\documentclass{article}
\usepackage{amsfonts}
\usepackage{amsmath}
\usepackage{amsthm}
\usepackage{url}
\usepackage{caption, wrapfig}
\usepackage{array}
\usepackage[usenames,dvipsnames]{xcolor}
\usepackage{todonotes}

\newtheorem{theorem}{Theorem}
\newtheorem{corollary}[theorem]{Corollary}
\newtheorem{proposition}[theorem]{Proposition}

\theoremstyle{definition}
\newtheorem{definition}[theorem]{Definition}

\newtheorem{example}[theorem]{Example}
\newtheorem{observation}[theorem]{Observation}

\usepackage[utf8]{inputenc}

\usepackage{tikz}
\usetikzlibrary{matrix, arrows, plotmarks}

\usepackage{etex}
\usepackage{multicol}
\usepackage[all]{xy}

\definecolor{fore}{RGB}{5,33,76}
\definecolor{back}{RGB}{245,222,179}
\definecolor{title}{RGB}{255,120,0}

\newcommand{\cmm}{\ensuremath{%
	\mathchoice
		{\includegraphics[height=1.7ex]{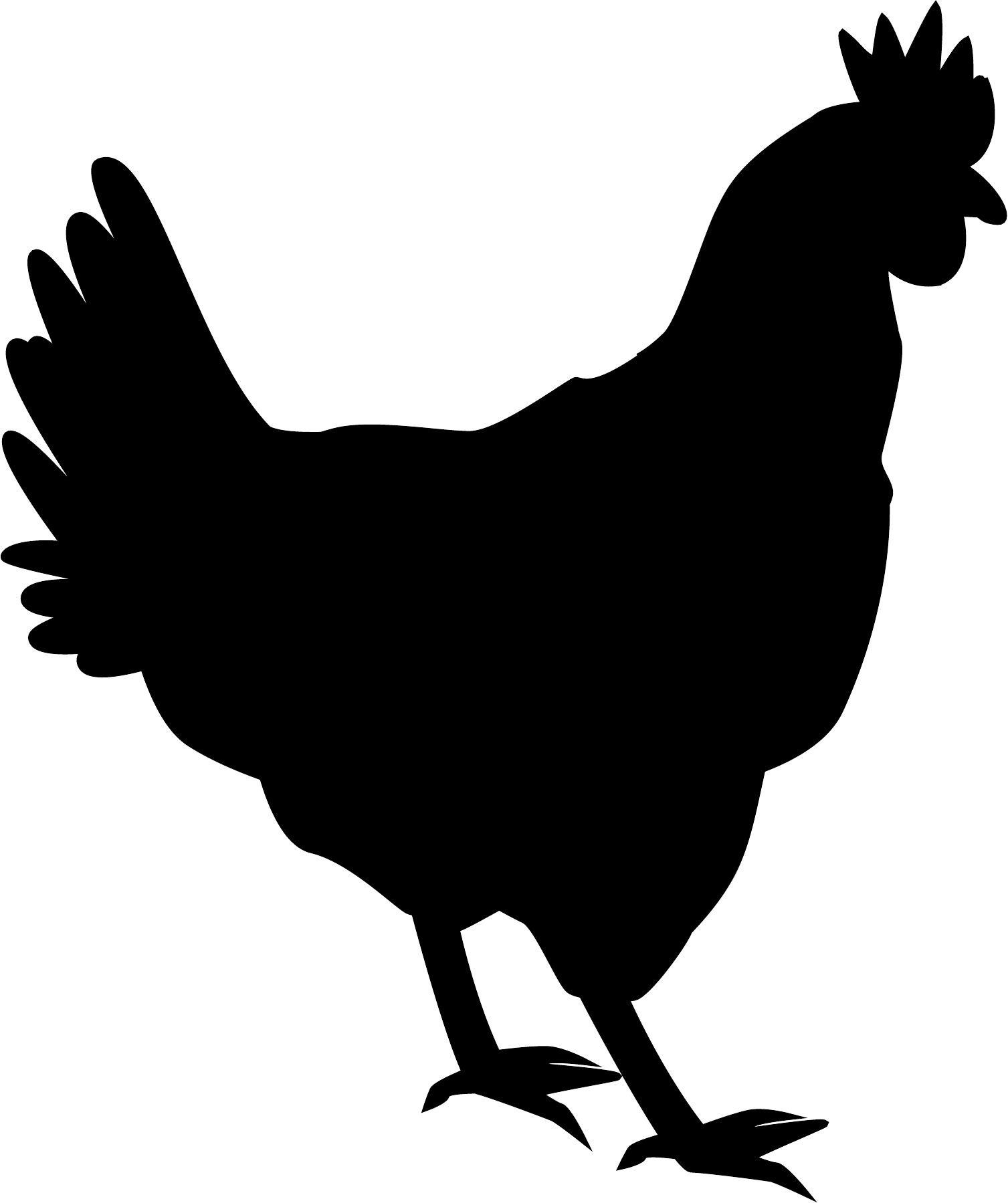}}
		{\includegraphics[height=1.7ex]{hen-311285.pdf}}
		{\includegraphics[height=1.5ex]{hen-311285.pdf}}
		{\includegraphics[height=1.0ex]{hen-311285.pdf}}
}}
\newcommand{\nm}{\langle n_1, \ldots, n_k\rangle}
\newcommand{\dis}{\mathbf{d}}

\title{Distances Between Factorizations in the Chicken McNugget Monoid}
\author{Scott Chapman, Pedro Garc\'ia-S\'anchez, and Christopher O'Neill}

\begin{document}

\maketitle

\begin{abstract}
We use the Chicken McNugget Monoid to demonstrate various factorization properties related to relations and chains of factorizations.  We study in depth the catenary and tame degrees of this monoid.
\end{abstract}

\bigskip

\begin{center}
\emph{Luck is a dividend of sweat. The more you sweat, the luckier you get.}\\ - Ray Kroc \cite{BQ}
\end{center}

\bigskip

%%%%%%%%%%%%%%%%%%%%%%%%%%%%%%%%%%%%%%%%%%%%%%%%%%%%%%%%%%%%%%%%%%%%%%%%%%%%%
\section{Prelude}%%%%%%%%%%%%%%%%%%%%%%%%%%%%%%%%%%%%%%%%%%%%%%%%%%%%%%%%%%%%
\label{sec:prelude}%%%%%%%%%%%%%%%%%%%%%%%%%%%%%%%%%%%%%%%%%%%%%%%%%%%%%%%%%%
%%%%%%%%%%%%%%%%%%%%%%%%%%%%%%%%%%%%%%%%%%%%%%%%%%%%%%%%%%%%%%%%%%%%%%%%%%%%%

In \cite{CON}, the authors examined in detail the \emph{Chicken McNugget Monoid} (denoted in that paper by $\cmm$) and its related factorization properties.  
\begin{wrapfigure}{l}{0\linewidth}
\includegraphics[width=2in]{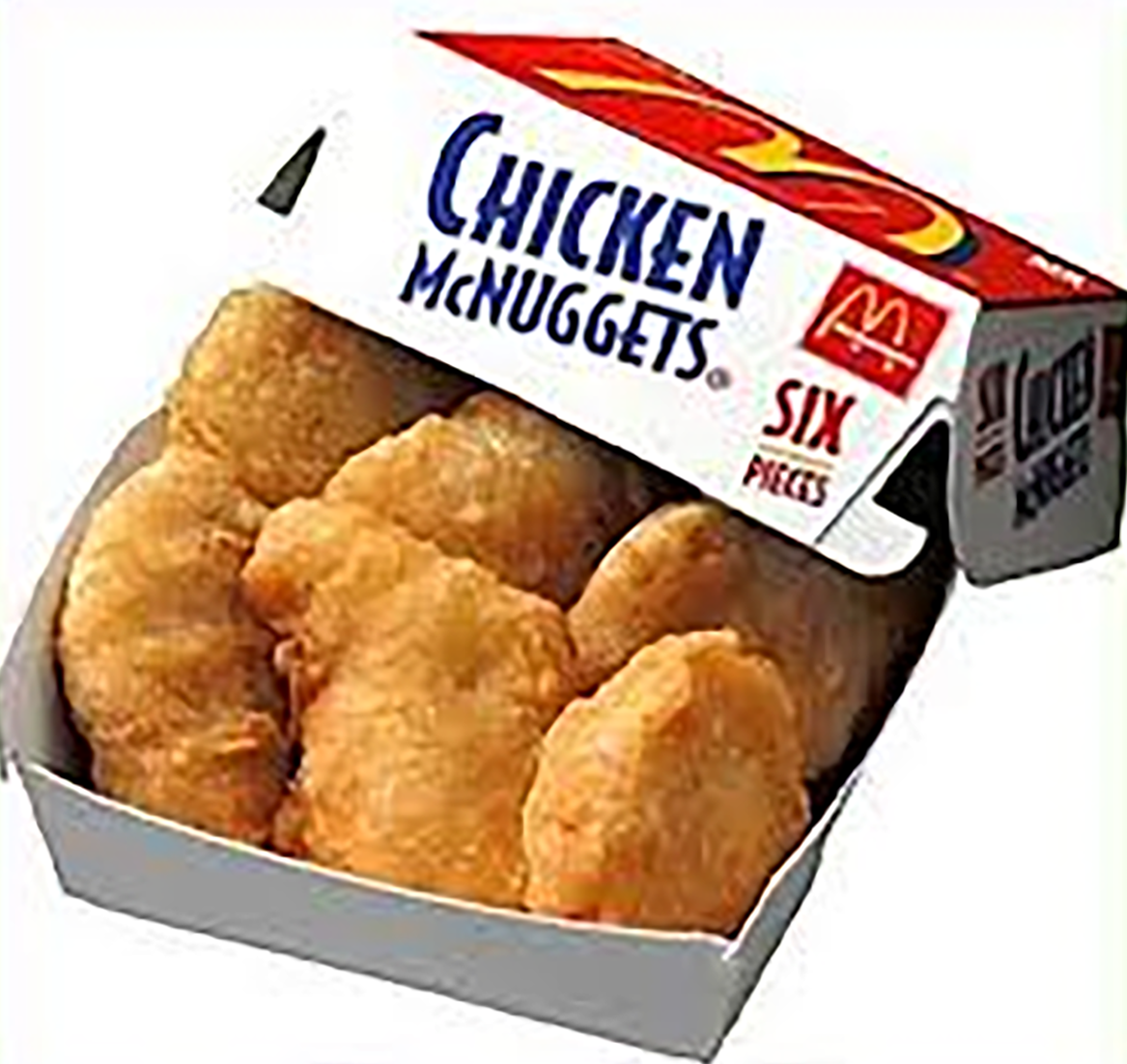}
\caption{The 6 piece box}
\end{wrapfigure}
These authors preceeded that paper with the following quote from McDonald's founder Ray Kroc \cite{BQ}: ``People just want more of it.''  From the reaction to that paper, Ray Kroc was right.  

In the present paper, we will pick up where \cite{CON} left off and explore \emph{chains of factorizations} in $\cmm$.  The notion of a chain of factorzations has up to now been largely relegated to research papers and not had wide exposure.  Our purpose here is to show, within the context of $\cmm$, that once the common factorization invariants such as elasticity, sets of length, and delta sets are determined, then a chain of factorizations, which relates to the complete set of factorizations of an element and not just the lengths, is a powerful tool in studying factorization properties.  We will introduce a method to measure the distance between two factorizations of a given element (see Definition \ref{d:distance}) and from this distance function compute two combinatorial constants:\ the \textit{catenary degree} (see Definition~\ref{d:catenary}) and the \textit{tame degree} (see Definition~\ref{d:tamedegree}).  While the catenary degree will in some sense measure the total ``spread'' of the distances between the complete set of factorizations of an element, the tame degree will focus on measuring distances from a factorization to another factorization containing a particular atom.  As with \cite{CON}, we present the definitions and examples in terms of a general numerical monoid, and conclude by specializing our results to the Chicken McNugget Monoid.

%%%%%%%%%%%%%%%%%%%%%%%%%%%%%%%%%%%%%%%%%%%%%%%%%%%%%%%%%%%%%%%%%%%%%%%%%%%%%
\section{Definitions and Basic Properties of the McNugget Monoid}%%%%%%%%%%%%
\label{sec:definitions}%%%%%%%%%%%%%%%%%%%%%%%%%%%%%%%%%%%%%%%%%%%%%%%%%%%%%%
%%%%%%%%%%%%%%%%%%%%%%%%%%%%%%%%%%%%%%%%%%%%%%%%%%%%%%%%%%%%%%%%%%%%%%%%%%%%%

So what is the Chicken McNugget Monoid?  We briefly review some background material which can be found in greater detail in \cite{CON}.  Chicken McNuggets were originally sold in packages of size 6, 9, or 20 pieces, and the question of how many Chicken McNuggets can be bought without breaking apart a package became a popular recreational mathematics question.  More specifically, if $n$ Chicken McNuggets can be purchased using whole packages (where $n$ is a positive integer), then there exist nonnegative integers $x_1, x_2, x_3 \in \mathbb{N}_0$ such that
\[
n = 6x_1 + 9x_2 + 20x_3.
\]
In this case, $n$ is called a \emph{McNugget number} and $(x_1, x_2, x_3)$ is called a \emph{McNugget expansion} of~$n$.  As $(30,0,0)$, $(0,20,0)$, $(15,10,0)$, and $(0,0,9)$ are all McNugget expansions of $n = 180$, it is clear that McNugget expansions of a given McNugget number need not be unique.  A full list of the McNugget expansions of McNugget numbers up to $n = 50$ can be found in \cite[Table~1]{CON}.  

Let 
\[
\langle 6,9,20\rangle = \{6x_1 + 9x_2 + 20x_3 : x_1, \, x_2, \, x_3\in \mathbb{N}_0\}
\]
represent the complete set of McNugget numbers.  Under regular integer addition, $\langle 6,9,20 \rangle$ forms a \emph{monoid}, meaning the sum of any two McNugget numbers is again a McNugget number.  As previously advertised, we will call this monoid the \emph{Chicken McNugget monoid} and denote it by $\cmm$.   In more generality, if $n_1,\ldots ,n_k$ is a set of relatively prime positive integers, then 
\[
\langle n_1, \ldots ,n_k\rangle = \{x_1n_1+x_2n_2+\cdots+x_kn_k\, : \,x_1, \ldots, x_k\in \mathbb{N}_0\}
\]
is known as a \emph{numerical monoid}.  A good general reference on numerical monoids (sometimes called numerical semigroups) is \cite{GSR}.  Given $n_1,\ldots ,n_k$ as above, it is easy using elementary number theory to argue that there is a largest positive integer not contained in $\langle n_1, \ldots ,n_k\rangle$.  This positive integer is known as the \emph{Frobenius number} of $\langle n_1, \ldots ,n_k\rangle$ and is the focus of much ongoing mathematics research (see \cite{RA}).  Using \cite[Proposition~1 and Table~1]{CON}, it follows that the Frobenius number of $\cmm$ is 43.  This is the largest number of Chicken McNuggets that cannot be ordered using whole boxes of sizes 6, 9, or 20.

\begin{wrapfigure}{l}{0\linewidth}
\includegraphics[width=2in]{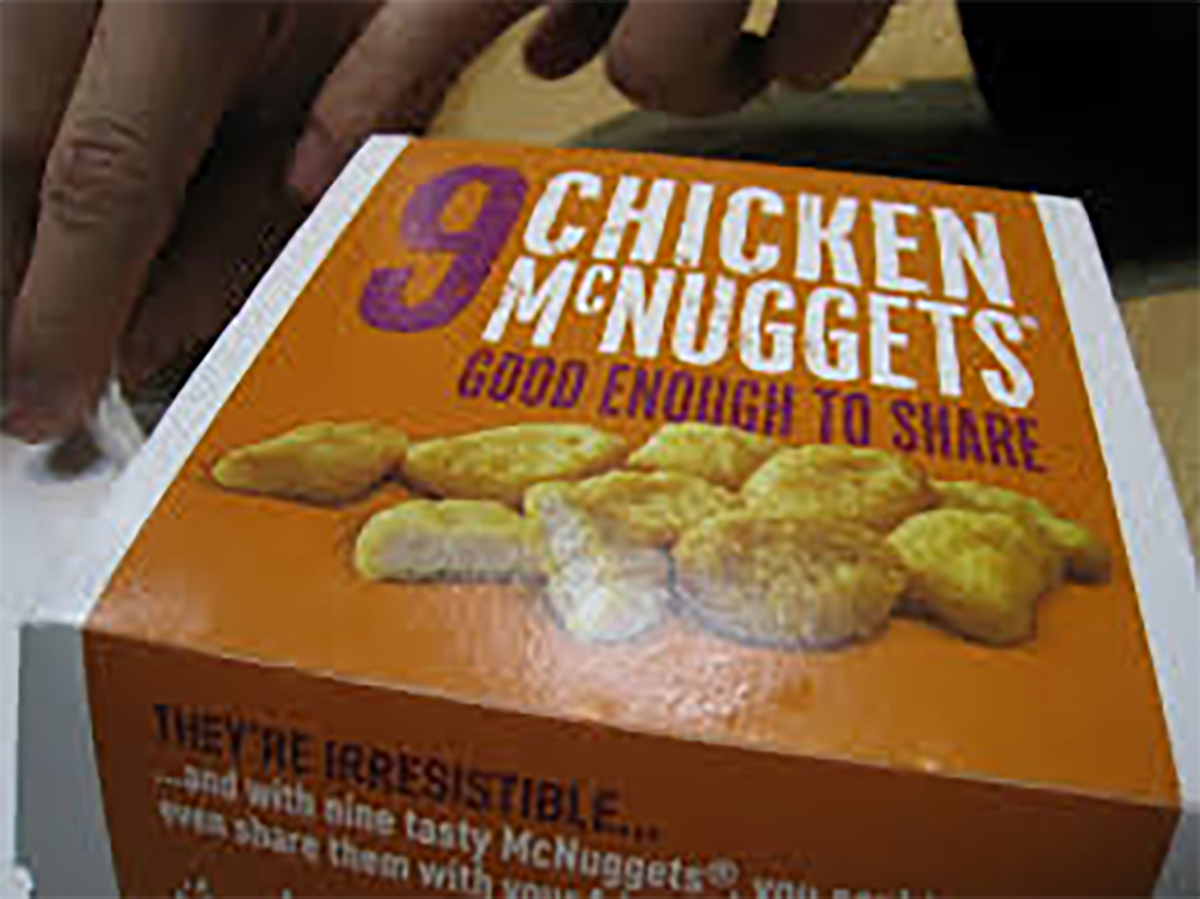}
\caption{The 9 piece box}
\end{wrapfigure}
In keeping with the usual notation used in commutative algebra, we will refer to the elements $n_1, \ldots ,n_k$ as the \emph{irreducible elements} or \emph{atoms} of $\nm$.  A representation of an element $n \in \nm$ as a sum $n = x_1n_1 + x_2n_2 + \cdots + x_kn_k$ of atoms will be called a \emph{factorization} of $n$ in $\nm$.  (Note that this is different from the ``usual'' notion of prime factorization of an integer!)  We will use the shorthand form $(x_1, \ldots, x_k)$ to represent a factorization of $n$ in $\nm$.  Set
\[
\mathsf{Z}(n) = \{ (x_1, \ldots, x_k) : n = x_1n_1 + x_2n_2 + \cdots + x_kn_k\}
\]
to be the complete set of factorizations of $n$ in $\nm$.  If $z = (x_1, \ldots, x_k) \in \mathsf{Z}(n)$, then the \textit{support} of $z$ is the set
$$\operatorname{supp}(z)=\{i\, :\, 1\leq i\leq k\mbox{  and  } x_i\neq 0\}.$$

Given a factorization $z = (x_1, \ldots, x_k) \in \mathsf{Z}(n)$, denote by $|z| = x_1 + \cdots + x_k$ the number of atoms (with repetition) used in $z$, called the \emph{length} of $z$.  The set
\[
\mathcal{L}(n)=\{\, |z|\, : \, z\in \mathsf{Z}(n)\}
\]
is known as the \emph{set of lengths} of $n$.  In the Chicken McNugget Monoid, each factorization $z$ of a McNugget number $n \in \cmm$ represents a specific combination of packs to purchase exactly $n$ McNuggets, and its length $|z|$ is simply the number of packs used.  

Most of the work in \cite{CON} centers around studying carefully defined \emph{invariants} that measure the size and structure of length sets of McNugget numbers.  
Writing the distinct lengths of a given element $n \in \cmm$ in order, we obtain $\mathcal{L}(n) = \{j_1, j_2, \ldots, j_t\}$ where $j_i < j_{i+1}$ for $i\in \{1,\ldots ,t-1\}$.  We further write 
\[
\ell(n) = j_1 \quad \text{and} \quad \operatorname{L}(n) = j_t
\]
for the minimum and maximum factorization lengths of $n$, respectively.  The \emph{elasticity} of $n$ is defined as the ratio
\[
\rho(n)=\frac{\operatorname{L}(n)}{\ell(n)},
\]
and the \emph{elasticity} of $\cmm$ as 
\[
\rho(\cmm)=\sup \{ \rho(n)\, : \, n\in \cmm\}.
\]
Intuitively, elasticity measures how ``spread out'' a monoid's factorization lengths~are.  The interested reader can find numerous papers in the recent literature that study problems related to elasticity, both in numerical monoids \cite{BOP1,BOP2,CHM} and more broadly \cite{DFA}.  

The \emph{delta set} of a McNugget number $n$ is defined by
\[
\Delta(n) = \{j_{i+1} - j_i : 1 \leq i \leq t-1\},
\]
and the \emph{delta set} of $\cmm$ by 
\[
\Delta(\cmm) = \bigcup_{n \in \cmm} \Delta(n).
\]
Intuitively, the delta set records the ``gaps'' in (or ``missing'') factorization lengths.  
There is a wealth of recent work concerning the computation of the delta set of a numerical monoid \cite{BOP2,BCKR,CDHK,CGLMS,CHK,CKLNZ,CK,GMV}.  For numerical monoids with three generators, the computation of the delta set is tightly related to Euclid's extended greatest common divisor algorithm \cite{GSLM1, GSLM2}.

We now summarize the main results in \cite[Corollary 9, Theorem 16]{CON}, which examine the elasticity and delta set of the Chicken McNugget Monoid.
\begin{proposition}\label{oldstuff}
Let $n\in \cmm$.
\begin{enumerate}
\item $\rho(\cmm ) = \frac{10}{3}$.
\item If $n \ge 92$, then
\[
\Delta(n) = \left\{\begin{array}{l@{\qquad}l}
\{1\} & \textnormal{if } r = 3, 8, 14, 17, \\
\{1, 2\} & \textnormal{if } r = 2, 5, 10, 11, 16, 19, \\
\{1, 3\} & \textnormal{if } r = 1, 4, 7, 12, 13, 18, \\
\{1, 4\} & \textnormal{if } r = 0, 6, 9, 15, \\
\end{array}\right.
\]
where $n = 20q + r$ for $q, r \in \mathbb{N}_0$ and $r < 20$.   
\item $\Delta(\cmm)=\{1,2,3,4\}$.  

\end{enumerate}
\end{proposition}

Before describing factorization chains, we note that the numerous calculations we will perform require some type of computing support.  The calculations we reference can be performed using the numericalsgps package \cite{DGM} for the computer algebra system GAP.  Interested readers are referred to that package for details behind the programming we use.

%%%%%%%%%%%%%%%%%%%%%%%%%%%%%%%%%%%%%%%%%%%%%%%%%%%%%%%%%%%%%%%%%%%%%%%%%%%%%
\section{Relations, Trades, and Minimal Presentations}%%%%%%%%%%%%%%%%%%%%%%%
\label{sec:minimalpresentations}%%%%%%%%%%%%%%%%%%%%%%%%%%%%%%%%%%%%%%%%%%%%%
%%%%%%%%%%%%%%%%%%%%%%%%%%%%%%%%%%%%%%%%%%%%%%%%%%%%%%%%%%%%%%%%%%%%%%%%%%%%%

We usually think of a numerical monoid $\nm$ in terms of its atoms $n_1, \ldots, n_k$.  Although these determine which integers live in $\nm$ and which do not, from the point of view of an algebraist, these only tell half of the story.  The underlying ``algebraic structure'' of $\nm$ also depends on the \emph{relations}, or linear dependencies, between $n_1, \ldots, n_k$.  

Let us examine what this means in the context of the Chicken McNugget Monoid.  The smallest McNugget number with more than one distinct McNugget expansion is~$18 \in \cmm$, which has $\mathsf{Z}(18) = \{(3,0,0), (0,2,0)\}$ since $18$ is a multiple of both $6$ and $9$.  This is precisely what is meant by a \emph{relation}, namely a linear equation relating the atoms $6$ and $9$.  This seemingly small observation has implications for most of the elements of $\cmm$;\ in \textbf{any} factorization $z = (x_1, x_2, x_3)$ of \textbf{any} McNugget number $n \in \cmm$, if $x_1 \ge 3$, then we can freely ``trade'' 3 copies of 6 for 2 copies of 9 to obtain another factorization of $n$, namely $(x_1 - 3, x_2 + 2, x_3)$.  We use the notation $(3,0,0) \sim (0,2,0)$ to represent this relation, indicating that in $\cmm$, $3$ times the first atom equals $2$ times the second.  

It is now natural to ask the following question.  Suppose you witness a customer ordering, say, $120$ Chicken McNuggets, using $10$ packs of $6$ and $3$ packs of $20$.  What~other ways are there to order that same number of Chicken McNuggets?  Well, using the relation $(3,0,0) \sim (0,2,0)$, we obtain at least 3 more ways, yielding
\[
(10,0,3), \, (7,2,3), \, (4,4,3), \, \text{and} \, (1,6,3).  
\]
Surely there must be others, since all of the above factorizations use the same number of $20$-packs.  To obtain these, we need another relation, one that involves the atom $20$.  Naturally, we should look for the smallest McNugget number that can be expressed using both packs of $20$ and packs of $6$ and/or $9$.  It turns out the magic number is $60 \in \cmm$, which has 
\[
\mathsf{Z}(60) = \{(10,0,0), (7,2,0), (4,4,0), (1,6,0), (0,0,3)\}.  
\]
We are now presented with a choice:\ which relation do we want?  Certainly it must involve the factorization $(0,0,3)$, but which of the $4$ factorizations involving $6$- and $9$-packs should be chosen?  Surprisingly, it does matter!  Whichever of the $4$ we choose will allow us to find all remaining factorizations of $120$ (and of any other McNugget number, for that matter).  

As an example, suppose we choose the relation $(10,0,0) \sim (0,0,3)$.  Starting with the initial factorization $(10,0,3)$, we can trade $6$'s for $20$'s to obtain $(0,0,6)$.  Moreover, we can instead trade $20$'s for $6$'s in $(10,0,3)$ and obtain 
\[
(20,0,0), \, (17,2,0), \, (14,4,0), \, (11,6,0), \, (8,8,0), \, (5,10,0), \, \text{and} \, (2,12,0)
\]
by subsequently trading $6$'s for $9$'s using $(3,0,0) \sim (0,2,0)$.  These turn out to be the final factorizations of $120 \in \cmm$.  Had we instead chosen to use the relation $(4,4,0) \sim (0,0,3)$, we can still obtain the factorization $(0,0,6)$, this time starting with $(4,4,3)$ and trading all of the $6$'s and $9$'s for $20$'s, and the remaining factorizations in the centered expression above can be obtained by swapping out the $20$'s in $(10,0,3)$, and then once again repeatedly applying $(3,0,0) \sim (0,2,0)$.  

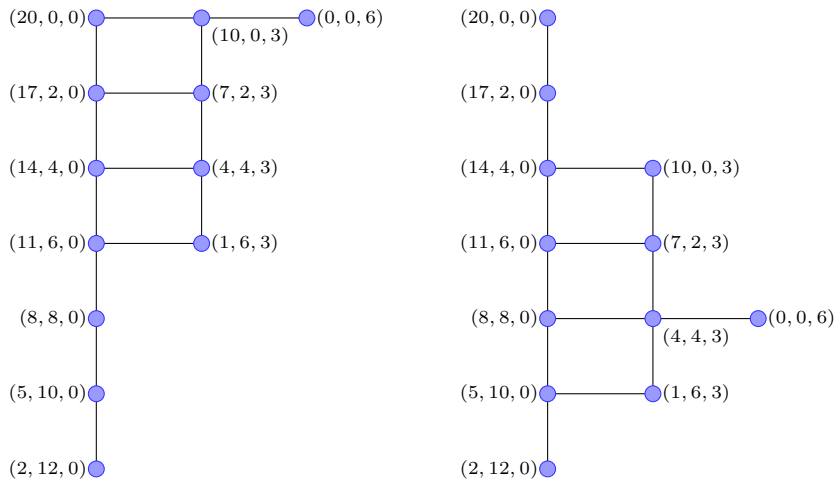
\begin{figure}
\begin{center}
\begin{tikzpicture}[y=.20cm, x=.20cm,font=\scriptsize]
\draw (0,0) to  (0,30);
\draw (0,30) to  (7,30);
\draw (0,25) to  (7,25);
\draw (0,20) to  (7,20);
\draw (0,15) to  (7,15) to (7,30) to (14,30);

\filldraw[fill=blue!40,draw=blue!80] (0,30) circle (3pt) node[anchor=east] {$(20,0,0)$};
\filldraw[fill=blue!40,draw=blue!80] (0,25) circle (3pt) node[anchor=east] {$(17,2,0)$};
\filldraw[fill=blue!40,draw=blue!80] (0,20) circle (3pt) node[anchor=east] {$(14,4,0)$};
\filldraw[fill=blue!40,draw=blue!80] (0,15) circle (3pt) node[anchor=east] {$(11,6,0)$};
\filldraw[fill=blue!40,draw=blue!80] (0,10) circle (3pt) node[anchor=east] {$(8,8,0)$};
\filldraw[fill=blue!40,draw=blue!80] (0,5) circle (3pt) node[anchor=east] {$(5,10,0)$};
\filldraw[fill=blue!40,draw=blue!80] (0,0) circle (3pt) node[anchor=east] {$(2,12,0)$};

\filldraw[fill=blue!40,draw=blue!80] (7,30) circle (3pt) node[anchor=north west] {$(10,0,3)$};
\filldraw[fill=blue!40,draw=blue!80] (7,25) circle (3pt) node[anchor=west] {$(7,2,3)$};
\filldraw[fill=blue!40,draw=blue!80] (7,20) circle (3pt) node[anchor=west] {$(4,4,3)$};
\filldraw[fill=blue!40,draw=blue!80] (7,15) circle (3pt) node[anchor=west] {$(1,6,3)$};

\filldraw[fill=blue!40,draw=blue!80] (14,30) circle (3pt) node[anchor=west] {$(0,0,6)$};

\draw (30,0) to  (30,30);
\draw (30,20) to  (37,20);
\draw (30,15) to  (37,15);
\draw (30,10) to  (37,10) to (44,10);
\draw (30,5) to  (37,5) to (37,20);

\filldraw[fill=blue!40,draw=blue!80] (30,30) circle (3pt) node[anchor=east] {$(20,0,0)$};
\filldraw[fill=blue!40,draw=blue!80] (30,25) circle (3pt) node[anchor=east] {$(17,2,0)$};
\filldraw[fill=blue!40,draw=blue!80] (30,20) circle (3pt) node[anchor=east] {$(14,4,0)$};
\filldraw[fill=blue!40,draw=blue!80] (30,15) circle (3pt) node[anchor=east] {$(11,6,0)$};
\filldraw[fill=blue!40,draw=blue!80] (30,10) circle (3pt) node[anchor=east] {$(8,8,0)$};
\filldraw[fill=blue!40,draw=blue!80] (30,5) circle (3pt) node[anchor=east] {$(5,10,0)$};
\filldraw[fill=blue!40,draw=blue!80] (30,0) circle (3pt) node[anchor=east] {$(2,12,0)$};

\filldraw[fill=blue!40,draw=blue!80] (37,20) circle (3pt) node[anchor=west] {$(10,0,3)$};
\filldraw[fill=blue!40,draw=blue!80] (37,15) circle (3pt) node[anchor=west] {$(7,2,3)$};
\filldraw[fill=blue!40,draw=blue!80] (37,10) circle (3pt) node[anchor=north west] {$(4,4,3)$};
\filldraw[fill=blue!40,draw=blue!80] (37,5) circle (3pt) node[anchor=west] {$(1,6,3)$};

\filldraw[fill=blue!40,draw=blue!80] (44,10) circle (3pt) node[anchor=west] {$(0,0,6)$};

\end{tikzpicture}
\end{center}
\caption{The factorization graph of $120$ in the Chicken McNugget Monoid with two different choices of minimal relations.}
\label{f:factgraphs120}
\end{figure}

It can be helpful to record the above information using a diagrams like the ones in Figure~\ref{f:factgraphs120} (called \emph{factorization graphs}).  Both graphs depict all of the factorizations of $120 \in \cmm$, but the left hand graph connects any two vertices with an edge if they are related by a single trade of $(3,0,0) \sim (0,2,0)$ or $(10,0,0) \sim (0,0,3)$, while the right hand graph uses the relations $(3,0,0) \sim (0,2,0)$ and $(4,4,0) \sim (0,0,3)$.  In both examples, we began at a factorization in the middle column, and used our second relation to branch out to the remaining columns.  Note that the edges of the factorization graph depend on a particular choice of relations.  

\begin{wrapfigure}{l}{0\linewidth}
\includegraphics[width=2in]{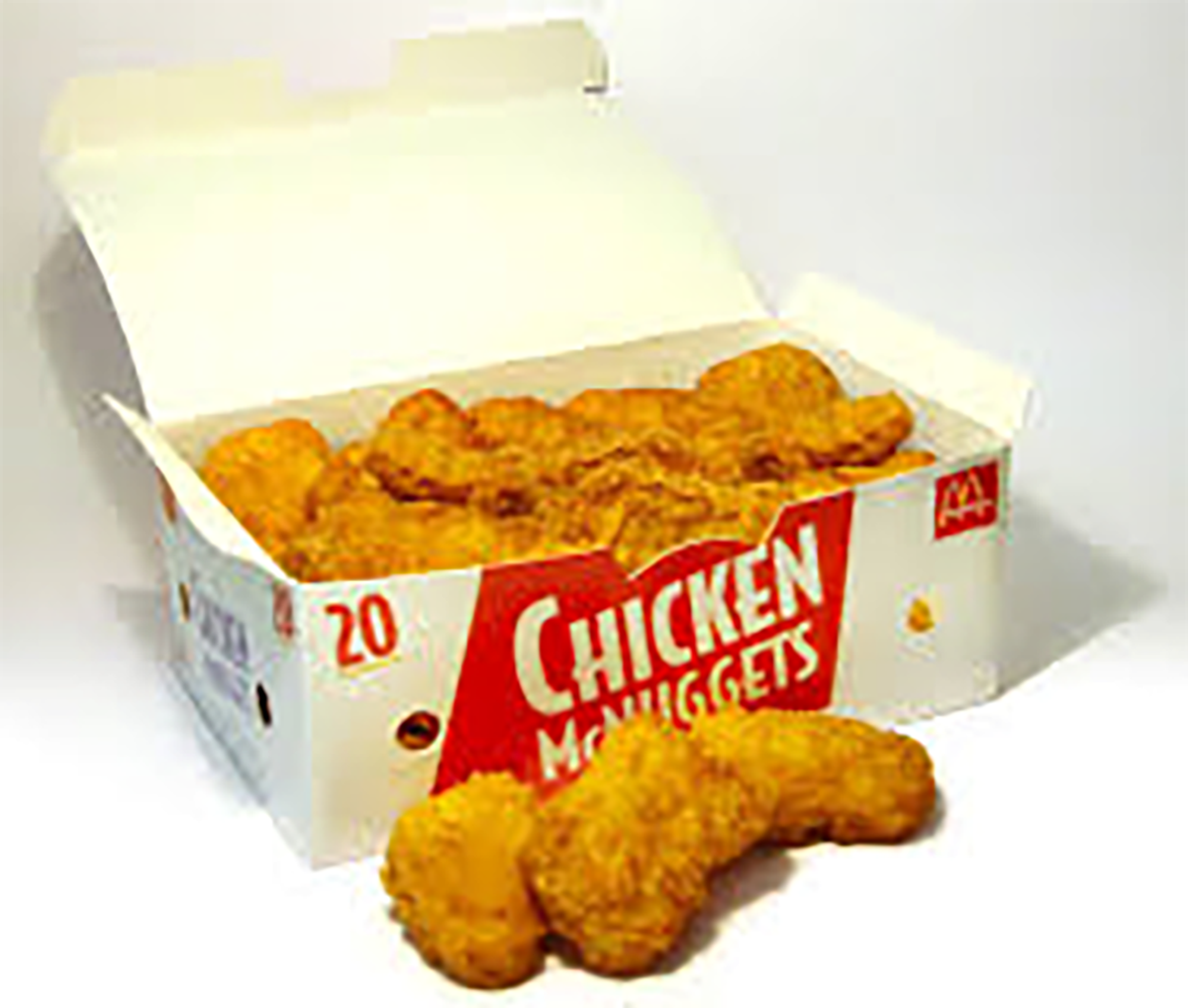}
\caption{The 20 piece box.}
\end{wrapfigure}

This illustrates the concept of a \emph{minimal presentation}, which is a collection $\rho$ of relations such that for an element $n$, any two factorizations of $n$ are connected by a sequence of trades using only the relations in $\rho$.  Said another way, a collection of relations forms a minimal presentation if for every~$n$, the factorization graph whose edges come from $\rho$ is connected.  The word ``minimal'' here means that none of the relations are implied by the rest (for instance, the relation $(6,0,0) \sim (0,4,0)$ would be redundant in $\cmm$ since it can be obtained by applying $(3,0,0) \sim (0,2,0)$ twice).  

Much is known about the structure of minimal presentations.  For instance, all of the minimal presentations for a given numerical monoid $\nm$ will have the same number of relations, and the elements whose factorizations appear in these relations will be identical as well.  Indeed, all $4$ minimal presentations for $\cmm$ involve one relation between factorizations of $18$ and one relation between factorizations of $60$.  One way to see this is that if we tried to build a minimal presentation $\rho$ using only the relation $(3,0,0) \sim (0,2,0)$, then $60$ would be the smallest element whose factorization graph was disconnected, implying that we must include in $\rho$ some relation between factorizations of $60$ to ensure its factorization graph is connected.  From there, as noted above, no matter which relation we pick, the factorization graphs of all remaining elements of $\cmm$ will be connected.  Indeed, this characterizes minimal presentations; they are minimal sets of relations so that any factorization graph is connected (see \cite[Chapter~7]{GSR} for thorough and precise definitions).  

Throughout the remainder of this paper,  we will use minimal presentations and factorization graphs to develop new invariants which will measure the relationships between the atoms of a general numerical monoid.  This will be completely analogous to how the elasiticy and delta set invariants measure the size and complexity of factorization lengths.  Along the way, we will encouter more graphs related to the factorization graph, but all will be different in key ways.  All of these graphs will be vital in our eventual arguments.

%%%%%%%%%%%%%%%%%%%%%%%%%%%%%%%%%%%%%%%%%%%%%%%%%%%%%%%%%%%%%%%%%%%%%%%%%%%%%
\section{The Amazing Distance Function}%%%%%%%%%%%%%%%%%%%%%%%%%%%%%%%%%%%%%%
\label{sec:distancefunction}%%%%%%%%%%%%%%%%%%%%%%%%%%%%%%%%%%%%%%%%%%%%%%%%%
%%%%%%%%%%%%%%%%%%%%%%%%%%%%%%%%%%%%%%%%%%%%%%%%%%%%%%%%%%%%%%%%%%%%%%%%%%%%%

In the previous section, we saw the role trades play in the structure of a numerical monoid $\nm$.  In order to define invariants from this structure, we need a way to measure which trades are ``larger'' than others.  Under such a measure, a trade $z_1 \sim z_2$ should measure as smaller than a trade that involves more atoms changing hands, but what does ``more'' mean?  

Given the important role that factorization lengths have played, it is tempting to consider the difference $\big||z_1| - |z_2|\big|$ in factorization lengths between $z_1$ and $z_2$ as a possible measure.  However, this has a drawback; consider the element $n = 126 \in \cmm$, which has factorizations 
\begin{multline*}
\mathsf{Z}(126) = \{(21,0,0), (18,2,0), (15,4,0), (12,6,0), (9,8,0), (6,10,0),\\ (3,12,0), (0,14,0), (11,0,3), (8,2,3), (5,4,3), (2,6,3), (1,0,6)\}.
\end{multline*}
Lurking in this set of factorizations is the trade $(11,0,3) \sim (0,14,0)$, which has a length difference of $0$, despite $14$ atoms being passed in each direction!  Clearly, this will not do.  

With this in mind, we consider the following measure of the ``size'' of a trade, one which focuses on the maximum length attained by the trade factorizations instead of their length difference.  We will make this definition in general terms so that it applies to all numerical moniods, and follow up with a concrete example.  

\begin{definition}\label{d:distance}
Let $S = \nm$ be a numerical monoid with $n \in S$.  Suppose $z_1 = (x_1,\ldots, x_k)$ and $z_2 = (y_1,\ldots ,y_k)$ are both in $\mathsf{Z}(n)$ and set 
\[
z_1 \wedge z_2 = (\min(x_1,y_1), \ldots, \min(x_n,y_n)).
\]
The \emph{distance} between the two factorizations $z_1$ and $z_2$ of $n$ is given by
\[
\dis(z_1, z_2) = \max\{|z_1|, |z_2|\} - |z_1 \wedge z_2|.
\]

\end{definition}
If $n = 126$, $z_1 = (0,14,0)$, $z_2 = (11,0,3)$, and $z_3 = (3,12,0)$, then we have
\[
\begin{array}{c}
\dis(z_1, z_2) = 14 - 0 = 14,\\
\dis(z_2, z_3) = 15 - 3 = 12, \\
\dis(z_1, z_3) = 15 - 12 = 3.
\end{array}
\]
Intuitively, $\dis(z_1, z_2)$ equals the maximum length of $z_1$ and $z_2$ where we have ignored the atoms appearing in both $z_1$ and $z_2$.  This ensures that a trade such as $(2,6,3) \sim (2,6,3)$ has distance $0$, which is reasonable considering that applying this trade has no net effect on the starting factorization.  

The distance function is an example of a \emph{metric}, meaning that it satisfies many of the same basic properties that other distances function do (you may have encountered metrics in an analysis class).  We gather some facts below and encourage the reader to work out their proofs as an exercise (the interested reader can also consult \cite[Proposition 1.2.5]{GHKb} for arguments).  

\begin{proposition}\label{p:distancefacts} 
If $S = \nm$ is a numerical monoid and $n \in S$ with $z_1, z_2, z_3 \in \mathsf{Z}(n)$, then we have the following:
\begin{enumerate}
\item $\dis (z_1,z_2)=0$ if and only if $z_1=z_2$;
\item if $z_1\neq z_2$, then $2\leq \dis (z_1,z_2)\leq \operatorname{L}(n)<\infty$;
\item $\dis (z_1, z_2)=\dis (z_2,z_1)$; and
\item $\dis (z_1, z_2)\leq \dis(z_1, z_3) + \dis (z_2,z_3)$.  
\end{enumerate}
\end{proposition}

The final item in Proposition~\ref{p:distancefacts} is known as the \emph{triangle inequality}, which, broadly speaking, ensures that one cannot find a strictly shorter distance between two points by first traveling to a third.  

We conclude this section with one last example, which will be used in the following section.  
In Table~\ref{tb:distexample}, we compute all the possible distances between factorizations of $103 \in \cmm$.  Due~to Proposition~\ref{p:distancefacts}.3, we need only fill in the top half of the table.
\begin{table}
\[
\arraycolsep=4.5pt\def\arraystretch{1.5}
\begin{array}{c|c|c|c|c}
& (9,1,2) & (6,3,2) & (3,5,2) & (0,7,2) \\
\hline
(9,1,2) &  0 & 3 & 6 & 9 \\
\hline
(6,3,2) &  & 0 & 3 & 6 \\
\hline
(3,5,2) &  & & 0 & 3 \\
\hline
(0,7,2) &  & & & 0  
\end{array}
\]
\caption{Distances between McNugget expansions of $103$.}
\label{tb:distexample}
\end{table}

%%%%%%%%%%%%%%%%%%%%%%%%%%%%%%%%%%%%%%%%%%%%%%%%%%%%%%%%%%%%%%%%%%%%%%%%%%%%%
\section{On Telephone Poles and Chains of Factorizations}%%%%%%%%%%%%%%%%%%%%
\label{sec:chains}%%%%%%%%%%%%%%%%%%%%%%%%%%%%%%%%%%%%%%%%%%%%%%%%%%%%%%%%%%%
%%%%%%%%%%%%%%%%%%%%%%%%%%%%%%%%%%%%%%%%%%%%%%%%%%%%%%%%%%%%%%%%%%%%%%%%%%%%%

So now that we know how to measure distances between factorizations, let us apply this to create an invariant which describes the distribution of the distances in $\mathsf{Z}(n)$.  

\begin{definition}\label{d:chain}
Let $S = \nm$ be a numerical monoid and $n \in S$.  A sequence of factorizations
\[
z_0, z_1, \ldots , z_t
\]
where each $z_i \in \mathsf{Z}(n)$ is called a \emph{chain of factorizations} of $n$.  For each $i\in\{1,\ldots,t\}$, set $\mathfrak{d}_i = \dis(z_{i-1}, z_i)$ which we refer to as the \emph{length} of the $i$-th link of the chain.
\end{definition}

Thus, in one sense you can think of a chain of the form $z_0, z_1, z_2, z_3, z_4$ in terms of the following picture, where the $\mathfrak{d}_i$'s represent the lengths of each individual ``link'' in the chain.  You can even think of the $\mathfrak{d}_i$'s as ``weights'' of the links.
\vskip12pt

\begin{tikzpicture}[y=.3cm, x=.35cm,font=\sffamily]

\draw[very thick, brown] (0,10) to[bend right] (8,10)  to[bend right] (16,10) to[bend right] (24,10) to[bend right] (32,10);

\filldraw[fill=blue!40,draw=blue!80] (0,10) circle (3pt) node[anchor=south] {$z_0$};
\filldraw[fill=blue!40,draw=blue!80] (8,10) circle (3pt) node[anchor=south] {$z_1$};
\filldraw[fill=blue!40,draw=blue!80] (16,10) circle (3pt) node[anchor=south] {$z_2$};
\filldraw[fill=blue!40,draw=blue!80] (24,10) circle (3pt) node[anchor=south] {$z_3$};
\filldraw[fill=blue!40,draw=blue!80] (32,10) circle (3pt) node[anchor=south] {$z_4$};

\node[below] at (4,10.5) {$\mathfrak{d}_1$};
\node[below] at (12,10.5) {$\mathfrak{d}_2$};
\node[below] at (20,10.5) {$\mathfrak{d}_3$};
\node[below] at (28,10.5) {$\mathfrak{d}_4$};

\end{tikzpicture}
\vskip12pt

\noindent
Given any two factorizations $z$ and $z^\prime$ of an element $n \in \nm$, one can build infinitely many chains between them, since in the definition of chain there is no stipulation that the $z_i$'s need be distinct.  We in some sense want to find a chain
linking $z$ and $z^\prime$ that uses links of minimal distance.  Hence, we introduce the following definition.  

\begin{definition}\label{d:nchain}
Let $S = \nm$ be a numerical monoid and $n \in S$, and let $N$ be a positive integer.  A chain of elements $z_0, z_1, \ldots, z_t$ in $\mathsf{Z}(n)$ is called an $N$-chain if each distance $\mathfrak{d}_i \le N$ for $i\in\{1,\ldots,t\}$.
\end{definition}

We extend the picture from above to provide examples of $N$-chains.  The chain image has been extended to a sequence of ``telephone poles'' labeled at the top with particular factorizations, and at the bottom by the trades being performed.  
We use factorizations of $103 \in\cmm$ from Table 1.  The following depicts a 9-chain from $(9,1,2)$ to $(0,7,2)$.
\vskip12pt

\begin{center}
\begin{tikzpicture}[y=.3cm, x=.33cm,font=\scriptsize]

\draw[very thick, brown] (0,0) -- (0,6);
\draw (0,6) to[bend right] (10,6);
\draw[very thick, brown] (10,0) -- (10,6);

\filldraw[fill=blue!40, draw=blue!80] (0,6) circle (3pt) node[anchor=south west] {$\!\!\!\!(9,1,2)$};
\filldraw[fill=blue!40, draw=blue!80] (10,6) circle (3pt) node[anchor=south] {$(0,7,2)$};

\filldraw[fill=blue!40, draw=blue!80] (0,0) circle (3pt);
\node[anchor=north] at (5,2) {$(9,0,0) \sim (0,6,0)$};
\filldraw[fill=blue!40, draw=blue!80] (10,0) circle (3pt);

\node[below] at (5,6.5) {$9$};

\end{tikzpicture}
\end{center}

\noindent
By substituting the link for one that passes through $(6,3,2)$, we obtain a 6-chain.  

\begin{center}
\begin{tikzpicture}[y=.3cm, x=.33cm,font=\scriptsize]

\draw[very thick, brown] (0,0) -- (0,6);
\draw (0,6) to[bend right] (10,6) to[bend right] (20,6);
\draw[very thick, brown] (10,0) -- (10,6);
\draw[very thick, brown] (20,0) -- (20,6);

\filldraw[fill=blue!40, draw=blue!80] (0,6) circle (3pt) node[anchor=south west] {$\!\!\!\!(9,1,2)$};
\filldraw[fill=blue!40, draw=blue!80] (10,6) circle (3pt) node[anchor=south] {$(6,3,2)$};
\filldraw[fill=blue!40, draw=blue!80] (20,6) circle (3pt) node[anchor=south] {$(0,7,2)$};

\filldraw[fill=blue!40, draw=blue!80] (0,0) circle (3pt);
\node[anchor=north] at (5,2) {$(3,0,0) \sim (0,2,0)$};
\filldraw[fill=blue!40, draw=blue!80] (10,0) circle (3pt);
\node[anchor=north] at (15,2) {$(6,0,0) \sim (0,4,0)$};
\filldraw[fill=blue!40, draw=blue!80] (20,0) circle (3pt);

\node [below] at (5,6.5) {$3$};
\node [below] at (15,6.5) {$6$};

\end{tikzpicture}
\end{center}

\noindent
By substituting one more time the link between $(6,3,2)$ and $(0,7,2)$ for one through $(3,5,2)$, we can reduce further to a 3-chain.

\begin{center}
\begin{tikzpicture}[y=.3cm, x=.33cm,font=\scriptsize]

\draw[very thick, brown] (0,0) -- (0,6);
\draw (0,6) to[bend right] (10,6) to[bend right] (20,6) to[bend right] (30,6);
\draw[very thick, brown] (10,0) -- (10,6);
\draw[very thick, brown] (20,0) -- (20,6);
\draw[very thick, brown] (30,0) -- (30,6);

\filldraw[fill=blue!40, draw=blue!80] (0,6) circle (3pt) node[anchor=south west] {$\!\!\!\!(9,1,2)$};
\filldraw[fill=blue!40, draw=blue!80] (10,6) circle (3pt) node[anchor=south] {$(6,3,2)$};
\filldraw[fill=blue!40, draw=blue!80] (20,6) circle (3pt) node[anchor=south] {$(3,5,2)$};
\filldraw[fill=blue!40, draw=blue!80] (30,6) circle (3pt) node[anchor=south] {$(0,7,2)$};

\filldraw[fill=blue!40, draw=blue!80] (0,0) circle (3pt);
\node[anchor=north] at (5,2) {$(3,0,0) \sim (0,2,0)$};
\filldraw[fill=blue!40, draw=blue!80] (10,0) circle (3pt);
\node[anchor=north] at (15,2) {$(3,0,0) \sim (0,2,0)$};
\filldraw[fill=blue!40, draw=blue!80] (20,0) circle (3pt);
\node[anchor=north] at (25,2) {$(3,0,0) \sim (0,2,0)$};
\filldraw[fill=blue!40, draw=blue!80] (30,0) circle (3pt);

\node [below] at (5,6.5) {$3$};
\node [below] at (15,6.5) {$3$};
\node [below] at (25,6.5) {$3$};
	
\end{tikzpicture}
\end{center}

\noindent
As such, even though $(9,1,2)$ and $(0,7,2)$ are distance $9$ apart, we can obtain one from the other using only trades with distance at most $3$.

%%%%%%%%%%%%%%%%%%%%%%%%%%%%%%%%%%%%%%%%%%%%%%%%%%%%%%%%%%%%%%%%%%%%%%%%%%%%%
\section{Distances Required to Build Chains:\ The Catenary Degree}%%%%%%%%%%%
\label{sec:catenarydegree}%%%%%%%%%%%%%%%%%%%%%%%%%%%%%%%%%%%%%%%%%%%%%%%%%%%
%%%%%%%%%%%%%%%%%%%%%%%%%%%%%%%%%%%%%%%%%%%%%%%%%%%%%%%%%%%%%%%%%%%%%%%%%%%%%

\begin{wrapfigure}{l}{0\linewidth}
\includegraphics[width=2in]{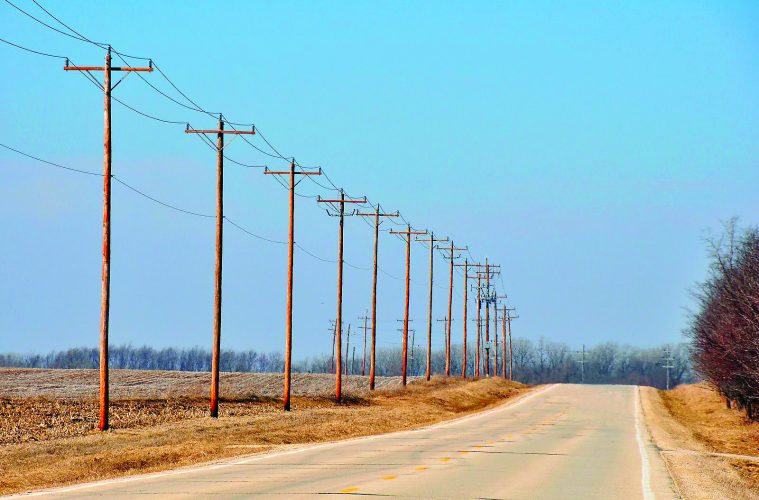}
\caption{Telephone poles never looked so good.}
\end{wrapfigure}

In the way that elasticity analyzes the ``spread'' of factorization lengths of an element, and the delta set analyzes the relative distribution within the set of lengths, how might one use the distance function to describe the structure of the set $\mathsf{Z}(n)$?  The answer lies in the $N$-chains constructed above.  In the case where factorizations in $\mathsf{Z}(n)$ are in close proximity to each other, one would expect to be able to construct an $N$-chain between any two factorizations for a small value of $N$.  The larger this necessary value of $N$, the more complex the structure of $\mathsf{Z}(n)$.  This motivates the following definition.

\begin{definition}\label{d:catenary}
Let $S = \nm$ be a numerical monoid and $n \in S$.  The \emph{catenary degree} of $n$ is defined as 
\[
\operatorname{c}(n) = \min\{ N : \mbox{there exists an } N\mbox{-chain beween any } z_1, z_2\in \mathsf{Z}(n) \}.
\]
We define the \emph{catenary degree} of $S$ to be
\[
\operatorname{c}(S) = \sup\{ \operatorname{c}(n) : n \in S \}.
\]
\end{definition}

\noindent 
Before continuing, we note that the study of the catenary degree in numerical monoids has been a frequent topic in the recent mathematical literature \cite{AGS,BGSG,CCMMP,CGSL,Om,OPTW}.  To understand some of the intricacies involved in studying the catenary degree, we will need to consider some of its elementary properties. 

Since the distance function cannot equal 1, $\operatorname{c}(n) = 0$ if and only if $|\mathsf{Z}(n)| = 1$ (that is, $n$ has unique factorization) and thus $\operatorname{c}(n) \ge 2$ if and only if $|\mathsf{Z}(n)| > 1$.  Moreover, it is easy to argue that $|\mathsf{Z}(n)| < \infty$ for every $n \in S$.  Thus the set 
\[
D(n) = \{ \dis(z_1, z_2) : z_1, z_2 \in \mathsf{Z}(n) \}
\]
is finite.  If $M > \max D(n)$, and $z_1$ and $z_2 \in \mathsf{Z}(n)$, then any chain from $z_1$ to $z_2$ is an $M$-chain, and hence $\operatorname{c}(n) < \infty$.  We summarize these fundamental observations in the next result.

\begin{proposition}\label{p:catenarybasic}
Let $S = \nm$ be a numerical monoid and $n \in S$.
\begin{enumerate}
\item $\operatorname{c}(n) = 0$ if and only if $|\mathsf{Z}(n)| = 1$.
\item If $|\mathsf{Z}(n)| \ne 1$, then $2 \le \operatorname{c}(n) < \infty$, and hence $\operatorname{c}(S) = 0$ or $2 \leq \operatorname{c}(S)$.
\end{enumerate}
\end{proposition}

It turns out that $\operatorname{c}(S)$ is always finite (and hence equal to the maximum of the catenary degrees achieved by the elements of $S$), though we defer a discussion on this matter until after the introduction of the tame degree in the next section.

While many of the references cited above work on computations of $\operatorname{c}(S)$, there is a relatively simple algorithm for obtaining $\operatorname{c}(n)$ from the set $\mathsf{Z}(n)$ using a graph similar to those used earlier.  Given $n \in S$, let $\mathcal{D}_n$ denote the complete graph whose vertices are the elements of $\mathsf{Z}(n)$, and label the edge between the factorizations $z_1$ and $z_2$ with $\dis(z_1,z_2)$.  We will refer to $\mathcal{D}_n$ as the \emph{distance graph} of $n$ with respect to $S$.  

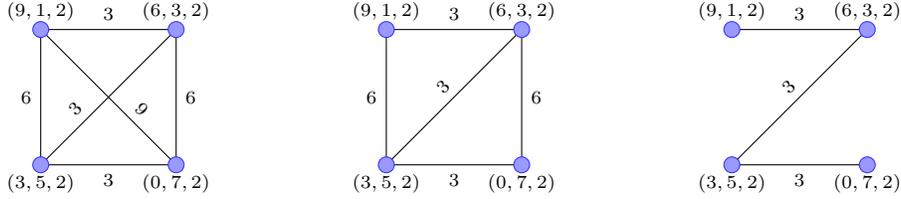
\begin{figure}
\begin{center}
\begin{tikzpicture}[y=.3cm, x=.3cm,font=\scriptsize]
\draw (0,0) -- (6,6);
\draw (0,0) -- (0,6);
\draw (0,0) -- (6,0);
\draw (6,0) -- (6,6);
\draw (6,0) -- (0,6);
\draw (0,6) -- (6,6);

\filldraw[fill=blue!40,draw=blue!80] (0,0) circle (3pt) node[anchor=north] {$(3,5,2)$};
\filldraw[fill=blue!40,draw=blue!80] (0,6) circle (3pt) node[anchor=south] {$(9,1,2)$};
\filldraw[fill=blue!40,draw=blue!80] (6,0) circle (3pt) node[anchor=north] {$(0,7,2)$};
\filldraw[fill=blue!40,draw=blue!80] (6,6) circle (3pt) node[anchor=south] {$(6,3,2)$};

\node [above] at (3,6) {$3$};
\node [below] at (3,0) {$3$};
\node [right] at (6,3) {$6$};
\node [left] at (0,3) {$6$};
\node [above, rotate=45] at (2,2) {$3$};
\node [above, rotate=-45] at (4,2) {$9$};
\end{tikzpicture}
\hfill
\begin{tikzpicture}[y=.3cm, x=.3cm,font=\scriptsize]
\draw (0,0) -- (6,6);
\draw (0,0) -- (0,6);
\draw (0,0) -- (6,0);
\draw (6,0) -- (6,6);
\draw (0,6) -- (6,6);

\filldraw[fill=blue!40,draw=blue!80] (0,0) circle (3pt) node[anchor=north] {$(3,5,2)$};
\filldraw[fill=blue!40,draw=blue!80] (0,6) circle (3pt) node[anchor=south] {$(9,1,2)$};
\filldraw[fill=blue!40,draw=blue!80] (6,0) circle (3pt) node[anchor=north] {$(0,7,2)$};
\filldraw[fill=blue!40,draw=blue!80] (6,6) circle (3pt) node[anchor=south] {$(6,3,2)$};

\node [above] at (3,6) {$3$};
\node [below] at (3,0) {$3$};
\node [right] at (6,3) {$6$};
\node [left] at (0,3) {$6$};
\node [above, rotate=45] at (3,3) {$3$};
\end{tikzpicture}
\hfill
\begin{tikzpicture}[y=.3cm, x=.3cm,font=\scriptsize]
\draw (0,0) -- (6,6);
\draw (0,0) -- (6,0);
\draw (0,6) -- (6,6);

\filldraw[fill=blue!40,draw=blue!80] (0,0) circle (3pt) node[anchor=north] {$(3,5,2)$};
\filldraw[fill=blue!40,draw=blue!80] (0,6) circle (3pt) node[anchor=south] {$(9,1,2)$};
\filldraw[fill=blue!40,draw=blue!80] (6,0) circle (3pt) node[anchor=north] {$(0,7,2)$};
\filldraw[fill=blue!40,draw=blue!80] (6,6) circle (3pt) node[anchor=south] {$(6,3,2)$};

\node [above] at (3,6) {$3$};
\node [below] at (3,0) {$3$};
\node [above, rotate=45] at (3,3) {$3$};
\end{tikzpicture}
\end{center}
\caption{The distance graph of $103$ in the Chicken McNugget monoid, in full (left) and with some edges removed (middle and right).  As removing any remaining edge would yield a disconnected graph, $\operatorname c(n) = 3$.}
\label{f:distgraph103}
\end{figure}

\begin{example}\label{e:distgraph}
Consider the distance graph of $103 \in \cmm$, depicted in Figure~\ref{f:distgraph103}.  One way to obtain the catenary degree is to remove edges from $D_{103}$, starting with those of highest weight, until removing a particular edge disconnects the graph (such edges are known as \emph{bridges}).  The weight of the last edge removed equals the catenary degree.  Several of the graphs resulting from this process are depicted alongside the full distance graph in Figure~\ref{f:distgraph103}.  Removing any one edge would disconnect the last graph, so the catenary degree of $103$ is $\operatorname{c}(103) = 3$.  Our implementation is essentially the well-known ``reverse-delete'' algorithm which first appeared in a paper by Kruskal \cite{Kr56}.
\end{example}

Before proceeding on to a definition of the tame degree, we outline in less formal language the meaning of the catenary degree.

\medskip
\noindent\fbox{\parbox{4.88in}{
Summary:\ Let $S = \nm$ be a numerical monoid and $n \in S$. 
\begin{enumerate}
\item $\operatorname{c}(n) = N$ means that $N$ is the smallest positive integer such that an $N$-chain exists between any two factorizations of $n$.
\item $\operatorname{c}(S) = N$ means that $N$ is the smallest positive integer such that given any element $m \in S$, an $N$-chain exists between any two factorizations of $m$.
\end{enumerate}
}}

%%%%%%%%%%%%%%%%%%%%%%%%%%%%%%%%%%%%%%%%%%%%%%%%%%%%%%%%%%%%%%%%%%%%%%%%%%%%%
\section{Distances Required to Reach Atoms: The Tame Degree}%%%%%%%%%%%%%%%%%%
\label{sec:tamedegree}%%%%%%%%%%%%%%%%%%%%%%%%%%%%%%%%%%%%%%%%%%%%%%%%%%%%%%%
%%%%%%%%%%%%%%%%%%%%%%%%%%%%%%%%%%%%%%%%%%%%%%%%%%%%%%%%%%%%%%%%%%%%%%%%%%%%%

While the catenary degree measures length in terms of chains, the tame degree measures distance from factorizations containing a specified atom.  
To motivate this invariant, we return to $126\in \cmm$ and note the set $\mathsf{Z}(126)$ computed earlier.  Notice that $(21,0,0)\in \mathsf{Z}(126)$ and does not contain any copies of the atom $20$.  How close is it to a factorization that does?  There are 5 such factorizations, and we list their distances from $(21,0,0)$ in the following diagram.

\begin{center}
\begin{tikzpicture}[y=.20cm, x=.20cm,font=\scriptsize]
\draw (0,0) to  (10,0);
\draw (10,0) to  (20,0);
\draw (10,0) to (0,10);
\draw (10,0) to (10,10); 
\draw (10,0) to (20,10);

\filldraw[fill=blue!40,draw=blue!80] (00,0) circle (3pt) node[anchor=north] {$(11,0,3)$};
\filldraw[fill=blue!40,draw=blue!80] (10,0) circle (3pt) node[anchor=north] {$(21,0,0)$};
\filldraw[fill=blue!40,draw=blue!80] (20,0) circle (3pt) node[anchor=north] {$(1,0,6)$};
\filldraw[fill=blue!40,draw=blue!80] (00,10) circle (3pt) node[anchor=south] {$(8,2,3)$};
\filldraw[fill=blue!40,draw=blue!80] (10,10) circle (3pt) node[anchor=south] {$(5,4,3)$};
\filldraw[fill=blue!40,draw=blue!80] (20,10) circle (3pt) node[anchor=south] {$(2,6,3)$};

\node [above] at (5,0) {$10$};
\node [above] at (15,0) {$20$};
\node [right] at (10,5) {$16$};
\node [above] at (5,6) {$13$};
\node [above] at (15,6) {$19$};
\end{tikzpicture}
\end{center}

\noindent So $(21,0,0)$ is at minimum $10$ units distance from any factorization of $126$ which contains a copy of $20$.  We invite the reader to repeat this process on the remaining $7$ factorizations of $126$ which do not contain a copy of $20$; you will find that each such factorization is 10 units (or less) away from a factorization with a copy of~$20$.  

Measuring minimal distances from an arbitrary factorization to one that contains a specific atom is the idea behind the tame degee.  We give the technical definition of the tame degree below.

\begin{definition}\label{d:tamedegree}
Let $S = \nm$ be a numerical monoid with $n \in S$.
\begin{enumerate}
\item For each $i$ with $n - n_i \in S$, denote by $\operatorname{t}(n, n_i)$ the minimum $t$ such that for every $z \in \mathsf \mathsf{Z}(n)$, there exists a factorization $z^\prime \in \mathsf{Z}(n)$ with $z^\prime = (y_1,\ldots ,y_n)$ where $y_i \neq 0$ and $\dis(z,z^\prime)\leq t$.  If $n - n_i \notin S$, then define $\operatorname{t}(n, n_i) = 0$.
\item The \emph{tame degree} of $n$ is $\operatorname{t}(n) = \max\{ \operatorname{t}(n, n_i) : 1 \leq i \leq k\}$.
\item The \emph{tame degree} of $S$ is $\operatorname{t}(S) = \sup\{ \operatorname{t}(n) : n \in S\}$.
\end{enumerate}
\end{definition}

Hence, to compute $\operatorname{t}(n,n_i)$, for every factorization in $\mathcal{D}_n$ where the $n_i$-th coordinate is zero, we compute the minimum distance to a factorization where that coordinate is nonzero.  Thus, returning to $126\in\cmm$, our previous work has shown that $\operatorname{t}(126, 20) = 10$.  Notice that this required 40 distance calculations.  In a similar fashion, we obtain $\operatorname{t}(126, 6) = 3$ and $\operatorname{t}(126, 9) = 7$, meaning 
\[
\operatorname{t}(126) = \max\{10, 3, 7\} = 10.
\]
How is one to interpret this?   Given any factorizaton $z\in \mathsf{Z}(126)$, you can ``tame'' (or ``keep apart'') any two factorizations of 126 containing an arbitrarily chosen atom with a whip of length 10.

We establish some elementary properties of the tame degree in the next proposition, and as earlier leave the proofs to the reader.

\begin{proposition}\label{p:tame}
Let $S = \nm$ be a numerical monoid and $n \in S$.
\begin{enumerate}
\item 
We have $\operatorname{t}(n) = 0$ if and only if all factorizations in $\mathsf{Z}(n)$ have identical support.

\item 
We have $\operatorname{t}(n) \leq \operatorname{L}(n) < \infty$.  Hence either $\operatorname{t}(n) = 0$ or $2 \leq \operatorname{t}(n) < \infty$.
\end{enumerate}
\end{proposition}

While we have shown above a simple algorithm using graphs to compute $\operatorname{c}(n)$ for $n\in S$, we note that the computation 
of $\operatorname{t}(n)$ is in general much more complicated and not as intuitive.  
Hence, we close this section with a summary of the various tame degree definitions in practical terms.

\medskip

\noindent\fbox{\parbox{4.88in}{
Summary:\ let $S = \nm$ be a numerical monoid and $n \in S$. 
\begin{enumerate}
\item 
$\operatorname{t}(n,n_i) = m$ means $m$ is the smallest nonnegative integer such that if $z \in \mathsf{Z}(n)$, then there is some factorization $z^\prime \in \mathsf{Z}(n)$ containing at least one copy of $n_i$ that is within distance $m$ of $z$.  

\item 
$\operatorname{t}(n) = m$ means $m$ is the smallest nonnegative integer such that if $z \in \mathsf{Z}(n)$, then for each $i\in\{1, \ldots, k\}$, there is some factorization $z^\prime \in \mathsf{Z}(n)$ containing at least one copy of $n_i$ that is within $m$ units of $z$.  

\item 
$\operatorname{t}(S) = m$ means $m$ is the smallest nonnegative integer such that if $n \in S$ and $z \in \mathsf{Z}(n)$, then for each $i \in\{1, \ldots, k\}$, there is some factorization $z^\prime \in \mathsf{Z}(n)$ containing at least one copy of $n_i$ that is within $m$ units of $z$.

\end{enumerate}
}}

%%%%%%%%%%%%%%%%%%%%%%%%%%%%%%%%%%%%%%%%%%%%%%%%%%%%%%%%%%%%%%%%%%%%%%%%%%%%%
\section{Computing Catenary and Tame Degrees of a Numerical Monoid}%%%%%%%%%%
\label{sec:finiteness}%%%%%%%%%%%%%%%%%%%%%%%%%%%%%%%%%%%%%%%%%%%%%%%%%%%%%%%
%%%%%%%%%%%%%%%%%%%%%%%%%%%%%%%%%%%%%%%%%%%%%%%%%%%%%%%%%%%%%%%%%%%%%%%%%%%%%

While we have argued in Propositions~\ref{p:catenarybasic} and~\ref{p:tame} that $\operatorname{c}(n)$ and $\operatorname{t}(n)$ are always finite, we have skirted the larger issue of the finiteness of $\operatorname{c}(S)$ and $\operatorname{t}(S)$.  To settle this point, we appeal to the following result proven by undergraduates in an NSF supported REU program from the summer of 2013.  

\begin{theorem}[{\cite[Theorem~3.1]{CCMMP}}]\label{t:periodic}
Let $S = \nm$ be a numerical monoid and suppose that $L = \mathrm{lcm}\{n_1, \ldots, n_k\}$.  The sequences $\{\operatorname{c}(n)\}_{n \in S}$ and $\{\operatorname{t}(n)\}_{n \in S}$ are eventually periodic with fundamental period a divisor of $L$.  
\end{theorem}

Thus, if $m$ is the point in $S$ at which $\{\operatorname{c}(n)\}_{n \in S}$ becomes periodic, then 
\[
\operatorname{c}(S) \in \{ \operatorname{c}(n) : n \in S, \, n \le m + L \}
\]
and hence must be finite.  Similar reasoning holds for $\operatorname{t}(S)$.  

\begin{corollary}\label{c:periodic}
If $S=\langle n_1, \ldots ,n_k\rangle$ is a numerical monoid, then both $\operatorname{c}(S)$ and $\operatorname{t}(S)$ are finite.
\end{corollary}

\begin{figure}[h]
\begin{center}
\includegraphics[width=5in]{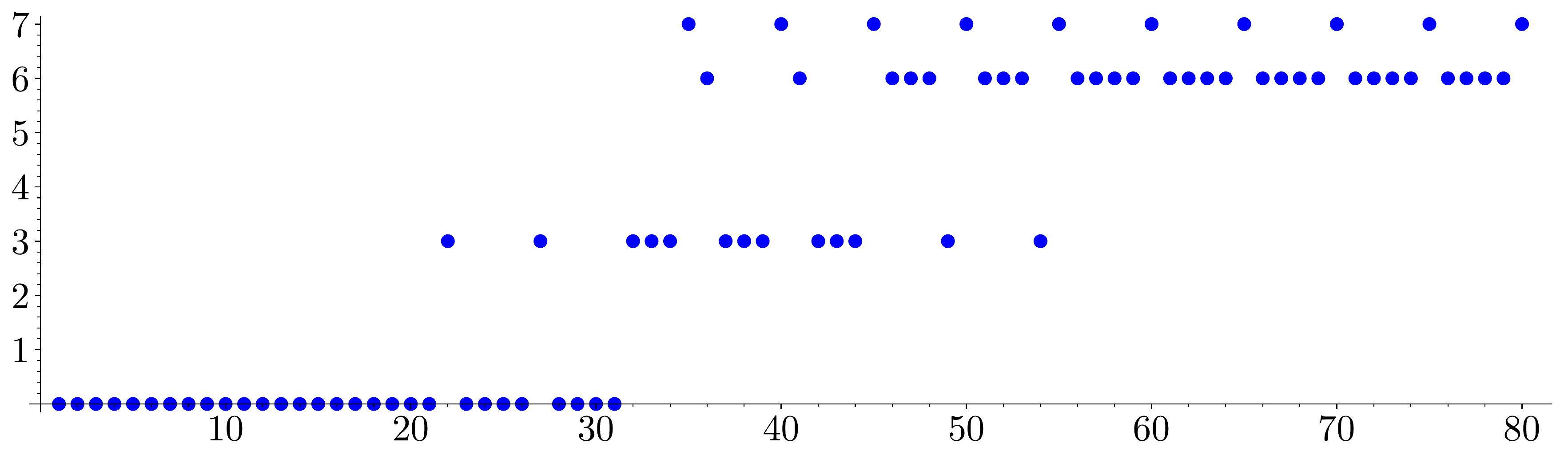}
\end{center}
\caption{A plot in which each point $(n, N)$ indicates $\operatorname{c}(n) = N$ for $n \in \langle 5, 11, 12 \rangle$.}
\label{f:periodic}
\end{figure}

\begin{example}\label{e:periodic}
The catenary degrees of the elements of $\langle 5, 11, 12 \rangle$ are depicted in Figure~\ref{f:periodic}.  One can readily observe that for $n \ge 55$, the catenary degree $\operatorname{c}(n)$ is periodic in $n$ with fundamental period $5$.  
\end{example}

Corollary~\ref{c:periodic} reduces the computation of $\operatorname{c}(S)$ and $\operatorname{t}(S)$ to a finite set of elements.  With a little more work we can do even better, restricting to so-called \emph{Betti elements} for the catenary degree and the \textit{Ap\'ery set} for the tame degree.  In the remainder of this section, we explore these constructions.

%%%%%%%%%%%%%%%%%%%%%%%%%%%%%%%%%%%%%%%%%%%%%%%%%%%%%%%%%%%%%%%%%%%%%%%%%%%%%
\subsection{Those beautiful Betti elements and awesome Ap\'ery sets}%%%%%%%%%
%%%%%%%%%%%%%%%%%%%%%%%%%%%%%%%%%%%%%%%%%%%%%%%%%%%%%%%%%%%%%%%%%%%%%%%%%%%%%

Let us return to the idea of \emph{minimal presentations} from earlier.  As we saw, given a numerical monoid $S = \nm$, a minimal presentation is a set of trades with which, for any $n \in S$, one can obtain any factorization in $\mathsf{Z}(n)$ from any other.  Using the language of chains, if $N$ is the highest trade distance in a minimal presentation of $S$, then there exists an $N$-chain between any two factorizations of $n$.  
% As such, if we want to compute the catenary degree, this $N$ is an upper bound.  
This allows us to identify which elements of $S$ are key to computing $\operatorname{c}(S)$.  

\begin{definition}\label{d:betti}
Let $S = \nm$ be a numerical monoid.  For $n \in S$, construct a graph $\mathcal{G}_n$, called the \emph{Betti graph}, whose vertices are the factorizations in $\mathsf{Z}(n)$, where an edge between $z_1$ and $z_2$ is included if $z_1$ and $z_2$ have at least one atom in common.  We call $n$ a \emph{Betti element} of $S$ if the Betti graph of $n$ is not connected.  Denote by $\mathrm{Betti}(S)$ the set of Betti elements of $S$.
\end{definition}

\begin{figure}[t]
\begin{center}
\begin{tikzpicture}[y=.20cm, x=.20cm,font=\scriptsize]
\filldraw[fill=blue!40,draw=blue!80] (0,7) circle (3pt) node[anchor=east] {$(3,0,0)$};
\filldraw[fill=blue!40,draw=blue!80] (0,2) circle (3pt) node[anchor=east] {$(0,2,0)$};

\draw (12,7) to  (12,2) to (17,0) to (17,9) to (12,7);
\draw (12,2) to (17,9);
\draw (12,7) to (17,0);

\filldraw[fill=blue!40,draw=blue!80] (12,7) circle (3pt) node[anchor=east] {$(7,2,0)$};
\filldraw[fill=blue!40,draw=blue!80] (12,2) circle (3pt) node[anchor=east] {$(4,4,0)$};
\filldraw[fill=blue!40,draw=blue!80] (17,0) circle (3pt) node[anchor=west] {$(1,6,0)$};
\filldraw[fill=blue!40,draw=blue!80] (17,9) circle (3pt) node[anchor=west] {$(10,0,0)$};
\filldraw[fill=blue!40,draw=blue!80] (22,5) circle (3pt) node[anchor=west] {$(0,0,3)$};

\draw (40,7) to (40,2) to (45,0) to (45,9) to (40,7) to (48.5,5);
\draw (48.5,5) to (40,2) to (45,9);
\draw (40,7) to (45,0) to (48.5,5) to (45,9);

\filldraw[fill=blue!40,draw=blue!80] (40,7) circle (3pt) node[anchor=east] {$(7,3,0)$};
\filldraw[fill=blue!40,draw=blue!80] (40,2) circle (3pt) node[anchor=east] {$(4,5,0)$};
\filldraw[fill=blue!40,draw=blue!80] (45,0) circle (3pt) node[anchor=west] {$(1,7,0)$};
\filldraw[fill=blue!40,draw=blue!80] (45,9) circle (3pt) node[anchor=west] {$(10,1,0)$};
\filldraw[fill=blue!40,draw=blue!80] (48.5,5) circle (3pt) node[anchor=west] {$(0,1,3)$};

\end{tikzpicture}
\end{center}
\caption{The Betti graphs of $18$ (left), $60$ (center), and $69$ (right) in the Chicken McNugget Monoid.}
\label{f:bettigraphscmm}
\end{figure}
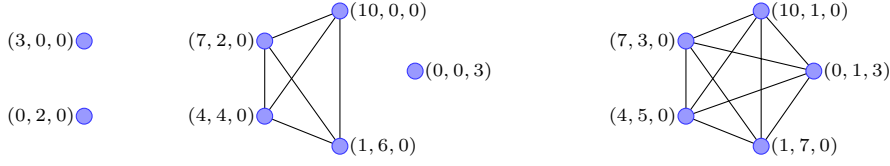

Returning to $\cmm$, we see in Figure~\ref{f:bettigraphscmm} that $\mathcal G_{18}$ consists of two vertices and no edges, and $\mathcal G_{60}$ has two connected components, one consisting of all factorizations involving $6$'s and $9$'s and the other a factorization using $20$'s.  Disconnected Betti graphs indicate that any minimal presentation must necessarily include a trade bridging the connected components.  On the other hand, Figure~\ref{f:bettigraphscmm} demonstrates that $\mathcal G_{69}$ is connected so $69$ is not a Betti element.  As it turns out, $\mathrm{Betti}(\cmm) = \{18,60\}$.  

In order to locate the Betti elements of $S = \nm$, we need to introduce a certain finite set of elements that sit at the heart of numerical monoids.  For motivation, consider the elements of $\cmm$ when organized based on their equivalence class modulo $6$:
\[
\cmm = \left\{
\begin{array}{@{}l@{\,\,}l@{\,\,}l@{\,\,}l@{\quad}l@{\,\,}l@{\,\,}l@{\,\,}l@{\quad}l@{\,\,}l@{\,\,}l@{\,\,}l@{}}
\textbf{0},  & 6,  & 12, & \ldots, &
\textbf{49}, & 55, & 61, & \ldots, &
\textbf{20}, & 26, & 32, & \ldots, \\
\textbf{9},  & 15, & 21, & \ldots, &
\textbf{40}, & 46, & 52, & \ldots, &
\textbf{29}, & 35, & 41, & \ldots
\end{array}
\right\}.
\]
Since $\cmm$ is closed under addition, every element of $\cmm$ can be obtained by adding a multiple of $6$ to one of the bolded values above, each of which is the smallest element of $\cmm$ in its equivalence class modulo $6$.  This leads to the following crucial definition.  

\begin{definition}
Let $S = \nm$ be a numerical monoid. For a nonzero $n \in S$, the \emph{Ap\'ery set} of $n$ in $S$, is defined and denoted as
\[
\operatorname{Ap}(S,n) = \{s \in S : s - n \not\in S\}.
\] 
\end{definition}

As discussed above, it is easy to see that there is a unique element in $\operatorname{Ap}(S,n)$ for each congruence class modulo $n$, each of which is precisely the minimum element of $S$ in its congruence class modulo $n$.  In particular, $|\operatorname{Ap}(S,n)| = n$.

Let us examine why Ap\'ery sets arise in the computation of the Betti elements. Assume that you have a bunch of factorizations, e.g., the factorizations
\[
\mathsf{Z}(60) = \{(10,0,0),(7,2,0),(4,4,0),(1,6,0),(0,0,3)\}
\]
of $60$.  In order to move from $(10,0,0)$ to $(7,2,0)$, we remove their ``common part'' $(7,0,0)$ and apply the trade $(3,0,0) \sim (0,2,0)$. Observe that $(3,0,0)$ and $(0,2,0)$ are factorizations of the same element $60 - (7 \cdot 6) = 18 \in \cmm$, and the factorizations
\[
\mathsf{Z}(18) = \{(3,0,0),(0,2,0)\}
\]
of $18$ have no common part, so there is no common part to remove. This means that the factorization where $6$ appears has no atom in common with any other factorizations (in this case, only $(0,2,0)$).  Hence $18 - 9$ is in $\cmm$, because we have a factorization where $9$ occurs, but $(18 - 9) - 6$ cannot be in $\cmm$, since this would imply that there is a factorization of $18$ where $6$ and $9$ both occur. This means $18 - 9 \in \operatorname{Ap}(\cmm,6)$, and $18 = n_i + w$ for $w = 18 - 9$ and $i \ne 1$.  

Notice that if we want to go from $(10,0,0)$ to $(4,4,0)$, the common part is $(4,0,0)$ and the new ``bridge'' is $6 \cdot 6 = 4 \cdot 9 = 36$, with factorizations 
\[
\mathsf{Z}(36) = \{(6,0,0),(3,2,0),(0,4,0)\}.
\]
Since we want to move from $(6,0,0)$ to $(0,4,0)$, we can use the fact that $(3,2,0)$ shares $3$ copies of $6$ with $(6,0,0)$ and $2$ copies of $9$ with $(0,4,0)$. In both situations, the problem reduces to moving from $(3,0,0)$ to $(0,2,0)$, which was already considered above as the factorizations of a Betti element. We can argue analogously with the rest of factorizations of $60$ that share some atoms, but there is one specific factorization, $(0,0,3)$, that does not share atoms with the rest. Since we can move freely now with trades in $\{(10,0,0),(7,2,0),(1,6,0),(4,4,0)\}$, it suffices to add a new trade to go from this set to $(0,0,3)$ (thus the different possible choices for minimal presentations for $\cmm$). Observe that in this case there is no factorization containing both $6$ and $20$. This means that $60-20\in \cmm$ but $(60-20)-6\not\in \cmm$, and as above $60=(60-20)+20$, with $60-20\in \operatorname{Ap}(\cmm,6)$ and $20$ a generator other than $6$.  This idea is behind the following result, which we will later find very useful.

\begin{theorem}\label{findbetti}\cite[Proposition~49]{AGSB}
Let $S=\nm$ be a numerical monoid minimally generated by $n_1,\ldots, n_k$ where $n_1<n_2<\cdots <n_k$.   If $s$ is a Betti element of $S$, then 
$s=n_i+w$ where $i\in \{2,\ldots ,k\}$ and $w\in \operatorname{Ap}(S, n_1)\backslash \{0\}$.
\end{theorem} 

We note that the converse of Theorem \ref{findbetti} is false; elements of the form $s=n_i+w$ in the theorem must be filtered before determining if they yield Betti elements.
By Theorem \ref{findbetti}, the computation of $\mathrm{Betti}(S)$ for a given $S = \nm$ is a finite process, but can be complicated, especially if $k$ is relatively large.  In fact, the size of $\mathrm{Betti}(S)$ can be arbitrarily large, even in the case $k = 4$ (in \cite{BR} a family with arbitrary number of Betti elements is given).  We will address this issue later, but for now we show why we are so interested in Betti elements.

\begin{theorem}\cite[Theorem 3.1]{CGSLPR}\label{t:betticatenary}
For any numerical monoid $S = \nm$, 
\[
\operatorname{c}(S) =\max\{\operatorname{c}(n) : n \in \mathrm{Betti}(S)\}.
\] 
\end{theorem}

\begin{example}\label{e:2gens}
We offer a very simple example to illustrate the ideas just presented.  Let $a$ and $b$ be relatively prime positive integers with $1 < a < b$, and set $S = \langle a, b \rangle$.  The elements of $S$ are of the form $ax + by$ where $x$ and $y$ are nonnegative integers.  Using Theorem \ref{findbetti}, it is easy to reason that
\[
\mathrm{Betti}(S) = \{ab\}\mbox{  and  }\operatorname{Ap}(S, a)=\{0,b,2b,\ldots, (a-1)b\}.
\]
Indeed, $\mathsf{Z}(ab) = \{(b,0), (0,a)\}$ and the Betti graph of any other element is either a single vertex (if $n - ab \notin S$) or connected (if $n - ab \in S$ is positive).  Thus, $\operatorname{c}(\langle a, b \rangle) = \operatorname{c}(ab)$.  Since $\dis((b,0), (0,a)) = b$, we conclude $\operatorname{c}(\langle a, b \rangle) = b$.
\end{example}

There is a somewhat similar method for computing $\operatorname{t}(S)$, though as with computing individual values of $\operatorname{t}(n)$, it is more expensive to complete.  
The method we will use centers around the following result. 

\begin{theorem}\cite[Theorem 16]{CGSL}\label{computetamecmm}
Let $S=\langle n_1,\ldots, n_k\rangle$ where the generating set for $S$ is minimal.  If $n$ is minimal in 
$S$ such that
$\operatorname{t}(n) = \operatorname{t}(S)$, then $n= w+n_i$ for some $i\in\{1,\ldots ,k\}$ and $w\in Ap(S,n_j)$ with $j\in\{1,\ldots ,k\}\backslash\{i\}$.
\end{theorem}

We note that there is an alternate method to compute $\operatorname{t}(S)$ which involves the computation of the \textit{primitive elements} of $\nm$, analogous to the Betti elements for the catenary degree; the interested reader should consult \cite[Proposition~4.1]{CGSLPR}.  

\begin{example}\label{twotame}
Returning to Example~\ref{e:2gens}, we again have that $\operatorname{t}(\langle a,b \rangle) = \operatorname{t}(ab)$, and as such, since $\dis((b,0),(0,a)) = b$, we conclude $\operatorname{t}(\langle a, b \rangle) = b$.  
\end{example}

We saw above that the Betti elements of $S$ were enough to compute the catenary degree of a numerical monoid, and these could be computed from the minimal generators and an Ap\'ery set. Thus computing the tame degree in general requires more machinery than computing the catenary degree.

%%%%%%%%%%%%%%%%%%%%%%%%%%%%%%%%%%%%%%%%%%%%%%%%%%%%%%%%%%%%%%%%%%%%%%%%%%%%%
\section{Calculations for the Chicken McNugget Monoid}%%%%%%%%%%%%%%%%%%%%%%%
\label{sec:computations}%%%%%%%%%%%%%%%%%%%%%%%%%%%%%%%%%%%%%%%%%%%%%%%%%%%%%
%%%%%%%%%%%%%%%%%%%%%%%%%%%%%%%%%%%%%%%%%%%%%%%%%%%%%%%%%%%%%%%%%%%%%%%%%%%%%

We begin with 
the Apéry set of $6\in \cmm$, which is 
\[
\operatorname{Ap}(S,6)=\{ 0, 49, 20, 9, 40, 29\},
\]  
written so the $i$-th element is the minimum element in $S$ congruent with $i$ modulo $6$. 
According to Theorem \ref{findbetti}, the candidates for Betti elements are 
\[
\{ 18, 29, 38, 40, 49, 58, 60, 69 \}.
\] 
We use GAP to find the factorizations of these elements, which are listed in Table~\ref{f:catenaryfacts}.  

\begin{table}[t]
\begin{tabular}{ll}
$n$ & $\mathsf{Z}(n)$ in $\cmm$\\ \hline
& \\[-5pt]
$18$ & $\{ ( 3, 0, 0 ), ( 0, 2, 0 ) \}$ \\
$29$ & $\{ ( 0,1,1 ) \}$ \\
$38$ & $\{ ( 3, 0, 1 ), ( 0, 2, 1 ) \}$ \\
$40$ & $\{ ( 0,0,2 ) \}$ \\
$49$ & $\{ ( 0,1,2 ) \}$ \\
$58$ & $\{ ( 3, 0, 2 ), ( 0, 2, 2 ) \}$ \\
$60$ & $\{ ( 10, 0, 0 ), ( 7, 2, 0 ), ( 4, 4, 0 ), ( 1, 6, 0 ), ( 0, 0, 3 ) \}$ \\
$69$ & $\{ ( 10, 1, 0 ), ( 7, 3, 0 ), ( 4, 5, 0 ), ( 1, 7, 0 ), ( 0, 1, 3 ) \}$ 
\end{tabular}
\caption{Factorizations of elements necessary to compute the catenary degree.}
\label{f:catenaryfacts}
\end{table}

In Figure~\ref{f:bettigraphscmm}, we have seen that $\mathcal{G}_{18}$ and $\mathcal{G}_{60}$ are disconnected and that $\mathcal{G}_{69}$ is connected.  The Betti graphs of 29, 40, and 49 are trivially connected as each is uniquely factorable, and those of 38 and 58 are connected by the trade $(3,0,0)\sim (0,2,0)$ alone.  Thus, the only Betti elements of $\cmm$ are 18 and 60.  The since $\mathcal{G}_{18}$ consists of two vertices, 
$\operatorname{c}(18) = \max\{2,3\} = 3$.    We can compute the catenary degree of 60 using the method outlined in Figure \ref{f:distgraph103}, and we reason through this proceedure using relations.  In order to move from any factorization to another in the set 
$\{(10,0,0),(7,2,0),(4,4,0),(1,6,0)\}$ we just need the relation $(3,0,0)\sim (0,2,0)$ which in terms of the catenary degree has a cost of three. And the shortest distance from this set to $(0,0,3)$ is attained by choosing $(1,6,0)$. This implies that $\operatorname{c}(60)=7$ and hence by Theorem~\ref{t:betticatenary}, we conclude $\operatorname c(S) = 7$.  

Now let us focus in the tame degree. According to Theorem~\ref{computetamecmm} we need to consider the factorizations of the elements in $n+\operatorname{Ap}(\cmm,m)$ for distinct $n, m \in \{6,9,20\}$. We already know $\operatorname{Ap}(\cmm,6)$; it is easy to check that 
\[
\operatorname{Ap}(\cmm, 10)=\{ 0, 46, 20, 12, 40, 32, 6, 52, 26 \},
\]
and 
\[
\operatorname{Ap}(\cmm, 20)=\{ 0, 21, 42, 63, 24, 45, 6, 27, 48, 9, 30, 51, 12, 33, 54, 15, 36, 57, 18, 39\}.
\]
So our set of elements of the form $n+w$ with $n$ a minimal generator of $\cmm$ and $w$ in the Ap\'ery set of another minimal generator is 
\begin{multline*}
\{  6, 9, 12, 15, 18, 20, 21, 24, 26, 27, 29, 30, 32, 33, 36, 38, 39, \\
40, 42, 45, 46, 48, 49, 51, 52, 54, 57, 58, 60, 63, 66, 69, 72 \}.
\end{multline*} 
Among these elements, $6, 9, 12, 15, 20, 21, 26, 29, 32, 40, 46, 49$, and $52$ each have a single factorization, and thus need not be considered.   The factorizations of the remaining elements can each be found in Table~\ref{f:catenaryfacts} or~\ref{tamefacts}.
\medskip

\begin{table}[t]
\begin{tabular}{ll}
$n$ & $\mathsf{Z}(n)$ in $\cmm$\\ \hline
& \\[-5pt]
$24$ & $\{ ( 4, 0, 0 ), ( 1, 2, 0 ) \}$ \\
$27$ & $\{ ( 3, 1, 0 ), ( 0, 3, 0 ) \}$ \\
$30$ & $\{ ( 5, 0, 0 ), ( 2, 2, 0 ) \}$ \\
$33$ & $\{ ( 4, 1, 0 ), ( 1, 3, 0 ) \}$ \\
$36$ & $\{ ( 6, 0, 0 ), ( 3, 2, 0 ), ( 0, 4, 0 ) \}$ \\
$39$ & $\{ ( 5, 1, 0 ), ( 2, 3, 0 ) \}$ \\
$42$ & $\{ ( 7, 0, 0 ), ( 4, 2, 0 ), ( 1, 4, 0 ) \}$ \\
$45$ & $\{ ( 6, 1, 0 ), ( 3, 3, 0 ), ( 0, 5, 0 ) \}$ \\
$48$ & $\{ ( 8, 0, 0 ), ( 5, 2, 0 ), ( 2, 4, 0 ) \}$ \\
$51$ & $\{ ( 7, 1, 0 ), ( 4, 3, 0 ), ( 1, 5, 0 ) \}$ \\
$54$ & $\{ ( 9, 0, 0 ), ( 6, 2, 0 ), ( 3, 4, 0 ), ( 0, 6, 0 ) \}$ \\
$57$ & $\{ ( 8, 1, 0 ), ( 5, 3, 0 ), ( 2, 5, 0 ) \}$ \\
$63$ & $\{ ( 9, 1, 0 ), ( 6, 3, 0 ), ( 3, 5, 0 ), ( 0, 7, 0 ) \}$ \\
$66$ & $\{ ( 11, 0, 0 ), ( 8, 2, 0 ), ( 5, 4, 0 ), ( 2, 6, 0 ), ( 1, 0, 3 ) \}$ \\
$69$ & $\{ ( 10, 1, 0 ), ( 7, 3, 0 ), ( 4, 5, 0 ), ( 1, 7, 0 ), ( 0, 1, 3 ) \}$ \\
$72$ & $\{ ( 12, 0, 0 ), ( 9, 2, 0 ), ( 6, 4, 0 ), ( 3, 6, 0 ), ( 0, 8, 0 ), ( 2, 0, 3 ) \}$
\end{tabular}
\caption{Factorizations of elements necessary to compute the tame degree.}
\label{tamefacts}
\end{table}

Observe that the tame degrees of $ 33$, $39$, $51$, and $57$ are each zero by Proposition~\ref{p:tame}, since all of their factorizations involve only the first two generators. The maximum distance between factorizations for $18$, $24$, $27$, $30$, $38$, $58$ is three. Notice that in the expressions of $36$, $42$, $45$, $48$, $54$, and $63$, only the first two generators appear (hence, these are acting like factorizations in the numerical monoid $\langle 2,3\rangle$ and the tame degree of this monoid is $3$; see Example~\ref{twotame}). Thus, the tame degrees of $18$, $24$, $27$, $30$, $36$, $38$, $42$, $45$, $48$, $54$, $58$, and $63$ are all $3$. 

So it remains to see what the tame degrees of $60$, $66$, $69$, and $72$ are.  We will only examine $60$ here, as the remaining elements can be handled in a similar fashion.  If we focus on the first generator, $6$, which appears in $(10,0,0)$, we have to find the closest factorization where $6$ does not occur.  The only candidate is $(0,0,3)$, and $\dis((10,0,0),(0,0,3)) = 10$. The distance between any other factorization where $6$ is involved and $(0,0,3)$ (the only one where $6$ does not occur) is less than $10$. But these factorizations are precisely those where $9$ appears, and so the tame degree does not grow when we look at the second generator. Now for the last generator, $20$, the only factorization in which it appears is $(0,0,3)$, and the closest where $20$ does not occur is $(1,6,0)$, and $\dis((0,0,3),(1,6,0)) = 7$. It follows that $\operatorname t(60) = 10$, and one can show that the same holds for $66$, $69$ and $72$. 

In total, we have obtained the following. 

\begin{proposition}\label{p:catenarytamemcnugget}
We have $c(\cmm) = 7$ and $t(\cmm) = 10$.
\end{proposition}

\begin{figure}[t]
\begin{center}
\includegraphics[width=5in]{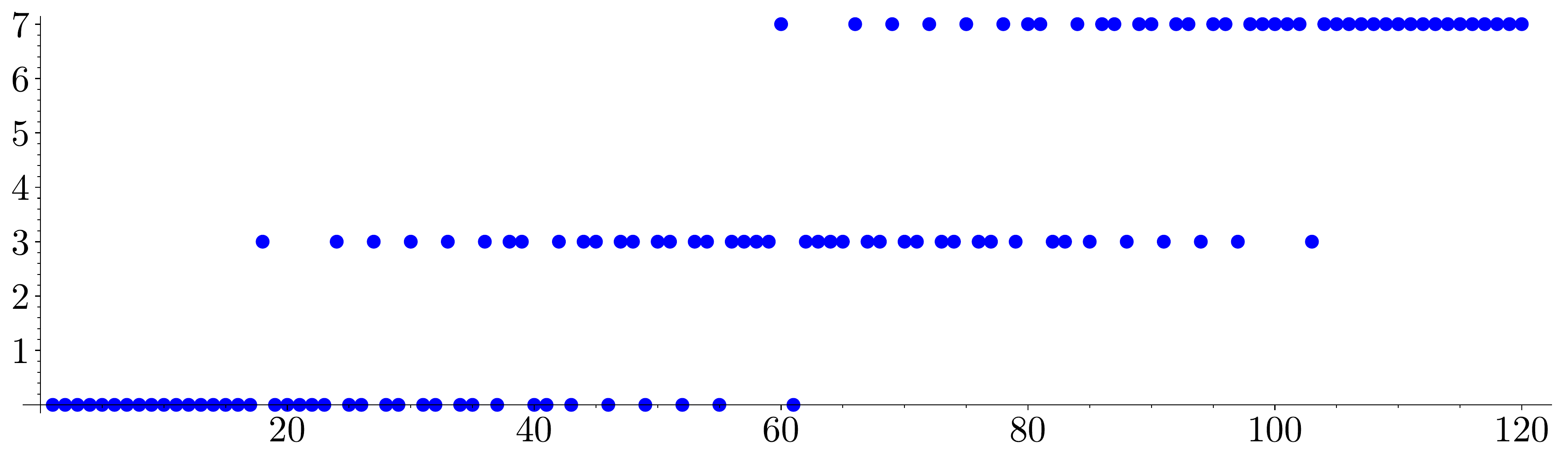}
\end{center}
\caption{A plot in which each point $(n, N)$ indicates $\operatorname{c}(n) = N$ for $n\in \langle 6,9,20\rangle$.}
\label{f:periodic1}
\end{figure}

\begin{figure}[t]
\begin{center}
\includegraphics[width=5in]{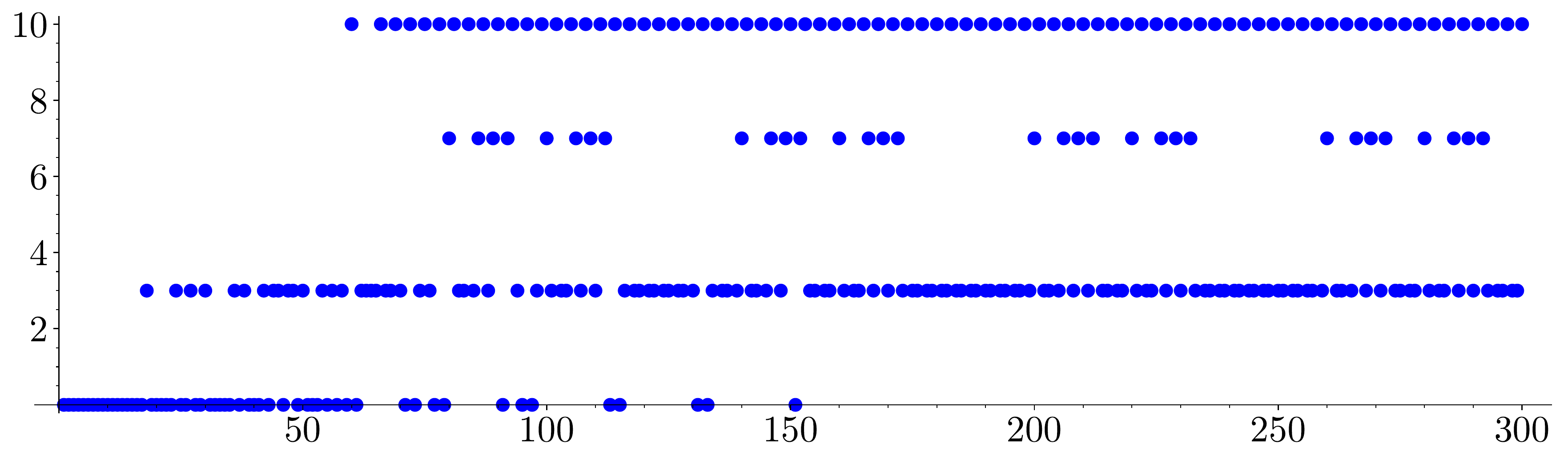}
\end{center}
\caption{A plot in which each point $(n, N)$ indicates $\operatorname{t}(n) = N$ for $n\in\langle 6,9,20\rangle$.}
\label{f:periodic2}
\end{figure}

We close by returning to Theorem~\ref{t:periodic} and give a complete description of the periodic behavior of the sequences $\{\operatorname{c}(s)\}_{s\in\cmm}$ and 
$\{\operatorname{t}(s)\}_{s\in\cmm}$.  Both sequences must have a fundamental period which divides $\operatorname{lcm}\{6,9,20\} = 180$, and Figures~\ref{f:periodic1} and~\ref{f:periodic2} give strong indication of the values indicated in Observation~\ref{o:periodicityspecifics}.  However, no proof of these observations are known aside from carefully examining factorizations and making arguments for each equivalence class modulo the fundamental periods, a particularly arduous task for the tame degree with its period of $60$.  

\medskip

\begin{observation}\label{o:periodicityspecifics}
The following hold.  
\begin{enumerate}
\item The sequence $\{\operatorname{c}(s)\}_{s\in\cmm}$ has fundamental period $1$ and begins at $n=104$.  Hence,
for $n\geq 104$, $\operatorname{c}(n)=7$.
\item The sequence  $\{\operatorname{t}(s)\}_{s\in\cmm}$ has fundamental period $60$ and begins at $n=152$.
\end{enumerate}
\end{observation}


\begin{thebibliography}{99}

\bibitem{AGS} F. Aguil\'o-Gost, Francesc, P. A. Garc\'ia-S\'anchez, Factorization and catenary degree in 3-generated numerical semigroups, \emph{Electron. Notes Discrete Math.} \textbf{34}(2009), 157-161.

\bibitem{DFA} D. F. Anderson, Elasticity of factorizations in integral domains: a survey, \textit{Lecture Notes in Pure Appl. Math.}, Marcel Dekker, New York, \textbf{189}(1997), 1--29.

\bibitem{AGSB} A. Assi, P. A.  Garc\'ia-S\'anchez, Numerical Semigroups and Applications,  RSME Springer series 1, Springer, Switzerland, 2016.

\bibitem{BOP1} T.~Barron, C.~O'Neill, R.~Pelayo, On the set of elasticities in numerical monoids, \emph{Semigroup Forum} \textbf{94} no.~1 (2017) 37--50.  
% Available at \textsf{arXiv:math.CO/1409.3425}.

\bibitem{BOP2}  T.~Barron, C.~O'Neill, R.~Pelayo, On dymamic algorithms for factorization invariants in numerical monoids, \emph{Math. Comp.} \textbf{86}(2017) 2429--2447.  
% Available at \textsf{arXiv:math.AC/1507.07435}

\bibitem{BGSG}  V. Blanco, P. A. Garcí\'a-S\'anchez, and A. Geroldinger, Semigroup-theoretical characterizations of arithmetical invariants with applications to numerical monoids and Krull monoids, \emph{Illinois J. Math.} \textbf{55}(2011) 1385-1414.

\bibitem{BCKR}  C. Bowles, S. T. Chapman, N. Kaplan, D. Reiser. On delta sets of numerical monoids. \emph{J. Algebra Appl.} \textbf{5}(2006) 695--718.

\bibitem{BR} H. Bresinsky, On prime ideals with generic zero $x_i = t^{n_i}$. \emph{Proc. Amer. Math. Soc.} \textbf{47} (1975) 329--332.

\bibitem{CCMMP} S. T. Chapman, M. Corrales, A. Miller, C. Miller, D. Patel, The catenary and tame degrees on a numerical monoid are eventually periodic, \emph{J. Aust. Math. Soc.} \textbf{97}(2014), 289--300. 

\bibitem{CDHK}  S. T. Chapman, J. Daigle, R. Hoyer, N. Kaplan,  Delta sets of numerical monoids using nonminimal sets of generators, \emph{Comm. Algebra} \textbf{38}(2010) 2622--2634.

\bibitem{CGSL} S. T. Chapman, P. A. Garc\'ia-S\'anchez, D. Llena, The catenary and tame degree of numerical monoids, \emph{Forum Math.} \textbf{21}(2009) 117--129.

\bibitem{CGSLPR} S. T. Chapman, P. A. Garc\'ia-S\'anchez, D. Llena, V. Ponomarenko, and J. C. Rosales, The catenary and
tame degree in finitely generated commutative cancellative monoids, \emph{Manuscripta Math.} \textbf{120}(3):253--
264, 2006.

\bibitem{CGLMS}  S. T. Chapman, P. A. Garc\'ia-S\'anchez, D. Llena, A. Malyshev, D. Steinberg, On the delta set and the Betti elements of a BF-monoid, \emph{Arabian J. of Mathematics} \textbf{1}(2012) 53--61.

\bibitem{CHM}  S. T. Chapman, M. T. Holden, T. A. Moore. Full elasticity in atomic monoids and integral domains, \emph{Rocky Mountain J. Math.} \textbf{37}(2006) 1437--1455.

\bibitem{CHK}  S. T. Chapman, R. Hoyer, N. Kaplan, Delta sets of numerical monoids are eventually periodic, \emph{Aequationes Math.} \textbf{77}(2009) 273--279.

\bibitem{CKLNZ}  S. T. Chapman, N. Kaplan, T. Lemburg, A. Niles, C. Zlogar, Shifts of generators and delta sets of numerical monoids, \emph{Inter. J. Algebra and Computation} \textbf{24}(2014) 655--669.

\bibitem{CON} S. T. Chapman and C. O'Neill, Factorization in the Chicken McNugget Monoid, \emph{Math. Magazine} \textbf{91}(2018), 323--336.

\bibitem{CK}  S. Colton, N. Kaplan, The realization problem for delta sets of numerical semigroups, \emph{J. Comm. Algebra}, to appear.

\bibitem{DGM}  M.~Delgado, P.~Garc\'ia-S\'anchez, J.~Morais, NumericalSgps, A package for numerical semigroups, Version 0.980 dev (2013), (GAP package), \url{http://www.fc.up.pt/cmup/mdelgado/numericalsgps/}.

\bibitem{Di} L. E. Dickson, Finiteness of the odd perfect and primitive abundant numbers with n distinct prime factors, \emph{Amer. J. Math.} 
\textbf{35}(4), 413–422.

\bibitem{Ga} J. Gallian, \emph{Contemporary Abstract Algebra}. Ninth Ed., Cengage Learning, Boston, MA, 2016.

\bibitem{GMV}  J. I. García-García, M. A. Moreno-Frías, A. Vigneron-Tenorio, Computation of Delta sets of numerical monoids, \emph{Monatshefte f\"ur Mathematik} \textbf{178}(2015) 457--472.

\bibitem{GSLM1} P. A. García-Sánchez, D. Llena, A. Moscariello, Delta sets for symmetric numerical semigroups with embedding dimension three, \emph{Aequationes Math.} 91 (2017), 579--600.

\bibitem{GSLM2} P. A. García-Sánchez, D. Llena, A. Moscariello, Delta sets for nonsymmetric numerical semigroups with embedding dimension three. \emph{Forum Math.} 30(1) (2018), 15--30.

\bibitem{GHKb} A.\ Geroldinger, F.\ Halter-Koch, \emph{Nonunique
    Factorizations: Algebraic, Combinatorial and Analytic Theory},
  Pure and Applied Mathematics, vol.\ 278, Chapman \& Hall/CRC, 2006.

\bibitem{BQ} Ray Kroc Quotes \url{https://www.brainyquote.com/authors/ray_kroc}

\bibitem{Kr56} J. Kruskal, On the shortest spanning subtree of a graph and the traveling salesman problem, \emph{Proc. Amer. Math. Soc.}, \textbf{7}(1956), 48--50. 

\bibitem{Om} M. Omidali, The catenary and tame degree of numerical monoids generated by generalized arithmetic sequences, \emph{Forum Math.} \textbf{24}(2012) 627--640.

\bibitem{OPTW} C. O’Neill, V. Ponomarenko, R. Tate, and G. Webb, On the set of catenary degrees of finitely generated cancellative commutative monoids, \textit{Internat. J.  Algebra Comput.}, \textbf{26}(2016), 565-576.

\bibitem{RA} J. L. Ram\'irez Alfons\'in, \emph{The Diophantine Frobenius Problem}, Oxford Lecture Series in Mathematics and Its Applications \textbf{30}, Oxford University Press, 2005, 256 pp.

\bibitem{GSR}  J. C. Rosales and P. A. Garc\'ia-S\'anchez, Numerical semigroups (Vol. 20). Springer Science \& Business Media, 2009.

\end{thebibliography}
\end{document}